\documentclass[oneside, reqno] {amsart}

\usepackage{xypic}
\input xy
\xyoption{all}
\usepackage{epsfig}
\usepackage{amsthm}
\usepackage{amssymb}
\usepackage{amsmath}
\usepackage{amscd}
\usepackage{color}
\usepackage[T1]{fontenc}
\usepackage{upgreek}
\usepackage[font=scriptsize]{caption}
\usepackage{wrapfig}
\usepackage{graphicx}
\usepackage{subfigure}

\newcommand{\bg}{\begin{equation}}
\newcommand{\ed}{\end{equation}}
\newcommand{\bga}{\begin{eqnarray}}
\newcommand{\eda}{\end{eqnarray}}
\newcommand{\pf}{\textbf{Proof:\ }}

\def\cbdu{\par{\raggedleft$\Box$\par}}

\newtheorem {Theorem}  {Theorem}

\numberwithin{Theorem}{section}

\newtheorem {Lemma}[Theorem]  {Lemma}

\theoremstyle{definition}

\theoremstyle{remark}
\newtheorem{Remark}[Theorem]{\bf Remark}

\expandafter\chardef\csname pre amssym.def
at\endcsname=\the\catcode`\@ \catcode`\@=11
\def\undefine#1{\let#1\undefined}
\def\newsymbol#1#2#3#4#5{\let\next@\relax
 \ifnum#2=\@ne\let\next@\msafam@\else
 \ifnum#2=\tw@\let\next@\msbfam@\fi\fi
 \mathchardef#1="#3\next@#4#5}
\def\mathhexbox@#1#2#3{\relax
 \ifmmode\mathpalette{}{\m@th\mathchar"#1#2#3}%
 \else\leavevmode\hbox{$\m@th\mathchar"#1#2#3$}\fi}
\def\hexnumber@#1{\ifcase#1 0\or 1\or 2\or 3\or 4\or 5\or 6\or 7\or 8\or
 9\or A\or B\or C\or D\or E\or F\fi}

\font\teneufm=eufm10 \font\seveneufm=eufm7 \font\fiveeufm=eufm5
\newfam\eufmfam
\textfont\eufmfam=\teneufm \scriptfont\eufmfam=\seveneufm
\scriptscriptfont\eufmfam=\fiveeufm

\catcode`\@=\csname pre amssym.def at\endcsname

\newcounter{remark}
\newcommand{\supp}{{\mathit supp}\,}

\newcommand{\R}{\mathbf{R}}

\def  \R   {{\mathbb R}}

\def  \12  {{\frac{1}{2}}}

\def\build#1_#2^#3{\mathrel{\mathop{\kern 0pt#1}\limits_{#2}^{#3}}}

\numberwithin{equation}{section}

\begin{document}
%\currannalsline{0}{2006}

\title[EMHD ill-posedness]{Ill-posedness of $2\frac12$D electron MHD}

%\author{hello}

\author [Mimi Dai]{Mimi Dai}

\address{Department of Mathematics, Statistics and Computer Science, University of Illinois at Chicago, Chicago, IL 60607, USA}
\email{mdai@uic.edu}

\thanks{The author is partially supported by the NSF grant DMS--2308208 and Simons Foundation.}

\begin{abstract}
We consider the electron magnetohydrodynamics (MHD) in the context where the 3D magnetic field depends only on the two horizontal plane variables. In particular, the magnetic field takes the form $B=\nabla\times (a\vec e_z)+b\vec e_z$ with $a=a(x,y)$ and $b=b(x,y)$.
Initial data $(a_0,b_0)$ is constructed in the Sobolev space $H^\beta \times H^{\beta-1}$ with $1<\beta<4$ such that the solution to this electron MHD system either escapes the space or develops norm inflation in $\dot H^\beta  \times \dot H^{\beta-1}$. 

\bigskip


\hspace{0.02cm}CLASSIFICATION CODE: 35Q35, 76E25, 76W05.
\end{abstract}

\maketitle

\section{Introduction}

%\medskip

%\subsection{Overview}

We study the electron magnetohydrodynamics (MHD) system
\begin{equation}\label{emhd}
\begin{split}
B_t+ \nabla\times ((\nabla\times B)\times B)=&\ 0,\\
\nabla\cdot B=&\ 0
\end{split}
\end{equation}
which is a simplified version of the full magnetohydrodynamics with Hall effect in the context of negligible ion flow motion and resistivity, see \cite{BDS, Bis1}. %The vector $B$ represents the magnetic field and $J=\nabla\times B$ the current density. The nonlinear term in (\ref{emhd}) captures the Hall effect which is believed to be responsible for the rapid magnetic reconnection phenomena in plasmas \cite{BDS}. 
%We note that $\nabla\cdot B(t)=0$ remains for all the time if $\nabla\cdot B(0)=0$ initially.
The unknown vector $B$ represents the magnetic field. 
Mathematical analysis of \eqref{emhd} is rather challenging albeit the model's importance in plasma physics. It is notorious that the nonlinear structure $\nabla\times ((\nabla\times B)\times B)$ in \eqref{emhd}, referred as the Hall term, is more singular than the nonlinear term $u\cdot \nabla u$ in the Euler equation which governs the pure hydrodynamics motion in the inviscid case. In particular, system \eqref{emhd} has the feature of being quasi-linear and supercritical (see \cite{Dai-W}).  

The investigation of well-posedness for \eqref{emhd} encounters great obstacles. Nevertheless well-posedness is addressed in the settings with a uniform magnetic background in \cite{Dai-emhd-2d, JO3}. On the other hand, ill-posedness phenomena have been discovered in different contexts in \cite{JO1, JO2}. The constructions in both \cite{JO1} and \cite{JO2} benefit from the dispersive structure of \eqref{emhd} and utilize degenerating wave packets techniques. In the current paper we pursue to establish ill-posedness for \eqref{emhd} in the two and half dimensional setting through a different approach. 

In physics literature, the two and half dimensional case of \eqref{emhd} is known to be of great importance. That is, 
\begin{equation}\label{two-half}
B(x,y,t)=\nabla\times (a\vec e_z)+b\vec e_z \ \ \mbox{with} \ \ \vec e_z=(0, 0,1), \ \ (x,y)\in\mathbb R^2,
\end{equation}
with scalar-valued functions
\[a=a(x,y,t), \ \ b=b(x,y,t).\]
Following \eqref{two-half}, it is clear that $B=(a_y, -a_x, b)$ and $\nabla\cdot B=0$. In term of $a$ and $b$, \eqref{emhd} can be written as the system
\begin{equation}\label{emhd-ab}
\begin{split}
a_t+ \nabla^{\perp}b\cdot \nabla a=&\ 0, \\
b_t+\nabla^{\perp}a\cdot\nabla \Delta a=&\ 0.
\end{split}
\end{equation}
We observe some interesting features of the system: 
\begin{itemize}
\item [(i)] the first equation \eqref{emhd-ab} is a transport equation, i.e. the horizontal potential is transported by the vorticity of the vertical component;
\item [(ii)] the nonlinearity in the second equation of \eqref{emhd-ab} is independent of $b$; 
\item [(iii)] the nonlinear term $\nabla^{\perp}a\cdot\nabla \Delta a=\{a,\Delta a\}$ is the Poisson bracket of $a$ and $\Delta a$; 
\item [(iv)] the stationary equation of the second one in \eqref{emhd-ab}
\[\nabla^{\perp}a\cdot\nabla \Delta a=0\]
coincides with the stationary Euler equation which has been studied to a great extent, see the recent work \cite{EHSX} and references therein.
\end{itemize}

The 2.5D electron MHD \eqref{emhd-ab} has been previously investigated by physicists, mostly using numerical simulations, for instance, see \cite{CSZ, KC} and references therein. The first mathematical study of \eqref{emhd-ab} with resistivity appeared in \cite{Dai-W}, where we showed the existence of determining modes and established a regularity criterion depending only on low modes of the solution. In particular, the low modes regularity condition implies the Beale-Kato-Majda (BKM) type blowup criterion is valid for the 2.5D electron MHD. It is known that the classical BKM blowup criterion is a vital property of the Euler equation and plays crucial roles in the study of singularity formation. In contrast, it remains an open question whether a BKM type blowup criterion can be established for \eqref{emhd} (or \eqref{two-half}) or not. Nevertheless, in the recent work \cite{Dai-Oh} we obtained a BKM type blowup criterion for the general 3D electron MHD \eqref{emhd} with resistivity $\Delta B$. 

In this paper we continue our investigation of the 2.5D electron MHD \eqref{emhd-ab} with the aim of establishing ill-posedness. 
The main result is stated below.

\begin{Theorem}\label{thm}
Let $1<\beta<4$. For arbitrarily small $\delta>0$, there exists an initial pair $(a_0,b_0)$ with $\|a_0\|_{H^\beta}+\|b_0\|_{H^{\beta-1}}\lesssim \delta$ such that the solution of \eqref{emhd-ab} with initial data $(a_0,b_0)$ develops norm inflation at a time $T\in(0,\delta]$, that is, 
\[\|a(T)\|_{\dot H^\beta}+\|b(T)\|_{\dot H^{\beta-1}}\gtrsim \frac1{\delta}.\]
\end{Theorem}

The symbols $\lesssim$ and $\gtrsim$ are explained in Subsection \ref{sec-notations} below.

\medskip

The crucial point in the construction of such norm inflation is to exploit the transport feature mentioned in (i). It is well-known that solutions of transport equation tend to lose regularity if the drift velocity is not Lipschitz, see the classical work \cite{DL, DM} for transport equation and Euler equation. Loss of smoothness of solutions to the 3D Euler equation due to non-Lipschitz velocity was further demonstrated in \cite{BT}. In the remarkable work \cite{BL1}, the authors established norm inflation for the Euler equation in borderline Sobolev spaces near Lipschitz space. 

In our case, we have a coupled system that consists a transport equation and an equation with highly singular nonlinear term. The main task is to construct initial data such that the component $a$ (which satisfies a transport equation) develops norm inflation, and in the same time to have the other component $b$ remain small. Inspired by the work \cite{CMZO} for the 2D Euler equation, we first seek a good approximating solution $(\bar a, \bar b)$ with $\bar a$ exhibiting norm inflation instantaneously; second, we perform perturbation analysis to show that $a-\bar a$ stays controlled up to the norm inflation time. 
To carry out the analysis for the coupled system, the features in items (ii), (iii) and (iv) also play important roles. 

It is worth to point out that the mechanism of establishing ill-posedness in the current paper is different from that of our previous work \cite{CD-continuity, CD-mhd, CD-norm} for the Navier-Stokes equations and the magnetohydrodynamics. In the latter cases, the main mechanism is the high-high-low interactions of frequencies which yields ill-posedness in spaces with negative derivative exponents.

\bigskip

\section{Preliminaries}
\label{sec-pre}

\subsection{Notations}
\label{sec-notations}
We use $C$ to denote a general constant which does not depend on other parameters in the text; it may be different from line to line.  The relaxed inequality symbol $\lesssim$ denotes $\leq$ up to a multiplication of constant when the constant is not necessary to be tracked. Analogously we use $\gtrsim$ as well.

\subsection{Equation in polar coordinate}
\label{sec-polar}
Our construction of initial data consists a radial part and an oscillation part. It is thus convenient to formulate the equation in the
standard polar coordinate $(r,\theta)$ with 
\[r=\sqrt{x^2+y^2}, \ \  x=r\cos\theta, \ \  y=r\sin \theta.\] 
Denote the polar coordinate unit basis vectors
\[e_r=(\cos\theta, \sin\theta), \ \ e_\theta=(-\sin\theta, \cos\theta).\]
For a function $f(x,y)$, 
%we also express it as \[f=f_r e_r+f_\theta e_\theta.\]
we have the following conversion of differentials
\begin{equation}\notag
\begin{split}
\partial_x f&=\cos\theta \ \partial_r f-\sin\theta\ \frac{\partial_\theta f}{r},\\
\partial_y f&=\sin\theta \ \partial_r f+\cos\theta\ \frac{\partial_\theta f}{r},\\
\nabla f&=\partial_r f e_r+\frac{\partial_\theta f}{r}e_\theta,\\
\nabla^\perp f&=\partial_r f\ e_\theta-\frac{\partial_\theta f}{r}\ e_r,\\
v\cdot\nabla f&=v_r\partial_r f+v_\theta\ \frac{\partial_\theta f}{r}.
\end{split}
\end{equation}

%For $f=f(\theta)$ independent of $r$, we have
%\[\Delta f=-\left[\frac12\sin(2\theta)+\frac1{r^2} \right]\partial_\theta f+\left[\cos^2\theta-\frac{\sin(2\theta)}{2r^2}\right]\partial_\theta^2f,\]
%\[\nabla^{\perp}f\cdot\nabla\Delta f=-\frac{2}{r^4}(\partial_\theta f)^2-\frac{\sin(2\theta)}{r^4}\partial_\theta f\partial_\theta^2f,\]
%\[\nabla^{\perp}(\nabla^{\perp}f\cdot\nabla\Delta f)=-\frac{2}{r^5}\left[...\right]\]

Thus a transport equation 
\[\partial_t f+u\cdot \nabla f=0\]
can be written in the polar coordinate form
\[\partial_t f+u_r\partial_r f+\frac{u_\theta}{r}\partial_\theta f=0\]
with $u=u_r e_r+u_\theta e_\theta$.

\bigskip

\section{Norm inflation of approximating solution}
\label{sec-approx}

Denote $u=\nabla^{\perp} b$. System \eqref{emhd-ab} can be written as 
\begin{equation}\label{emhd-au}
\begin{split}
\partial_t a+u\cdot \nabla a=&\ 0,\\
\partial_t u+\nabla^{\perp}(\nabla^{\perp}a\cdot\nabla\Delta a)=&\ 0.
\end{split}
\end{equation}

We will construct initial data $(a_0,b_0)$ such that $a_0$ consists a radial part and an oscillation part, while $b_0$ has only radial part. In that case, $u_0=\nabla^\perp b_0$ has only angular part. We then consider the approximating solution $(\bar a, \bar u)$ satisfying
\begin{equation}\label{app-a}
\begin{split}
\partial_t \bar a+u_0\cdot\nabla \bar a=&\ 0,\\
\bar a(x,0)=&\ a_0,
\end{split}
\end{equation}
and 
\begin{equation}\label{app-u}
\begin{split}
\partial_t \bar u+\nabla^{\perp}(\nabla^{\perp}\bar a\cdot\nabla\Delta\bar a)=&\ 0,\\
\bar u(x,0)=&\ u_0.
\end{split}
\end{equation}
The aim is to choose $(a_0,b_0)$ such that $\bar a$ develops norm inflation immediately in $H^\beta$; while $\bar u$ is left with the freedom to either remain controlled in $H^{\beta-2}$ for a short time or also develop norm inflation. 

Recalling the equation of $a$ in \eqref{emhd-ab}, we have the basic energy law
\begin{equation}\notag%\label{energy}
\frac12\frac{d}{dt}\int \left(a_x^2+a_y^2+b^2\right)\, dxdy=0
\end{equation}
%with $\vec x=(x,y)$
which indicates the a priori energy estimates 
\begin{equation}\notag%\label{basic-space}
 a\in L^\infty(0,T; H^1), \ \ 
 b\in L^\infty(0,T; L^2).
\end{equation}
Thus we pursue to show norm inflation of $(a,b)$ in $\dot H^\beta\times \dot H^{\beta-1}$ for $\beta>1$.
%According to the natural scaling of (\ref{emhd-ab}), the critical Sobolev space for $(a,b)$ is $\dot H^2\times \dot H^1$.
%\[a_\lambda=\lambda^{-1} a(\lambda \vec x, \lambda^2t), \ \ b_\lambda= b(\lambda \vec x, \lambda^2t).\]
%In two spatial dimension, some critical spaces for $a$ are 
%\[\dot H^2\subset \dot B^{1+\frac2p}_{p, \infty}\subset \dot B^{1}_{\infty, \infty}, \ \ 2\leq p<\infty, \]
%and some critical spaces for $b$ are
%\[\dot H^1\subset \dot B^{\frac2p}_{p, \infty}\subset \dot B^{0}_{\infty, \infty}, \ \ 2\leq p<\infty.\]
%Thus system (\ref{emhd-ab}) is energy supercritical even it is in 2D. The question of global existence of regular solution for (\ref{emhd-ab}) in general setting remains open. In \cite{Dai-W} the authors established a regularity criterion which only relies on the low modes of the solution. 

\medskip

\subsection{Initial data}
\label{sec-initial}
Let $\lambda\gg1$ be a large parameter, $g(r)$ and $h(r)$ be radial functions satisfying $g, h\in C_c^\infty (1, 4)$ and 
\[h'= 1 \ \ \mbox{on} \ \ (2,3).\] 
Consider the initial data
\begin{equation}\label{data}
a_0(r,\theta)=\delta\lambda^{1-\beta}g(\lambda r)\cos(\theta),\ \ b_0(r,\theta)=\delta\lambda^{2-\beta}h(\lambda r)
\end{equation}
for parameters $1<\beta<4$ and sufficiently small $\delta> 0$.
% {\color{blue} The radial part $\lambda^{1-\beta}f(\lambda r)$ may be not necessary.}
It follows that 
\begin{equation}\label{initial-u}
u_0(r,\theta)=\nabla^{\perp}b_0(r,\theta)=(\partial_rb_0) e_{\theta}=\delta\lambda^{3-\beta} h'(\lambda r) e_{\theta}.
\end{equation}

\begin{Lemma}\label{le-initial}
The estimates 
\begin{equation}\label{est-u0}
\|a_0\|_{H^s}\lesssim \delta\lambda^{s-\beta}, \quad
\|b_0\|_{H^s}\lesssim \delta\lambda^{s+1-\beta}, \quad
\|u_0\|_{H^s}\lesssim \delta\lambda^{s+2-\beta}
\end{equation}
hold for $s\geq 0$. In particular, we have
\[\|a_0\|_{H^\beta}+ \|b_0\|_{H^{\beta-1}}\lesssim \delta.\]
%\[a_0\in H^\beta, \ \ b_0\in H^{\beta-1}, \ \ u_0\in H^{\beta-2}.\]
Moreover, $u_0$ satisfies
\[\|u_0\|_{C^1}\approx \delta\lambda^{4-\beta}.\]

\end{Lemma}
\pf
It is straightforward to compute
%\begin{equation}\notag
%\begin{split}
%\int_0^{2\pi}\int_0^\infty \lambda^{2(1-\beta)}f^2(\lambda r)r\, drd\theta&= \lambda^{-2\beta}\int_0^{2\pi}\int_0^\infty f^2(\lambda r)(\lambda r)\, d(\lambda r)d\theta\\
%&\lesssim \lambda^{-2\beta}
%\end{split}
%\end{equation}
%and 
\begin{equation}\notag
\begin{split}
&\int_0^{2\pi}\int_0^\infty \lambda^{2-2\beta}g^2(\lambda r)\cos^2(\theta)r\, drd\theta\\
\leq&\ \lambda^{-2\beta}\int_0^{2\pi}\int_0^\infty g^2(\lambda r)(\lambda r)\, d(\lambda r)d\theta\\
\lesssim&\ \lambda^{-2\beta}
\end{split}
\end{equation}
which implies
\[ \|a_0\|_{L^2}\lesssim \delta\lambda^{-\beta}.\]
To estimate the higher norm $H^s$ of $a_0$,  we define
\[
F(x,y):=g(r)\frac{x}{r},\qquad (x, y)\in\mathbb R^2.
\]
Since $g\in C_c^\infty(1,4)$, we have $F\in C_c^\infty(\mathbb R^2)$ and
\[
a_0(x,y)=\delta\lambda^{1-\beta}F(\lambda x, \lambda y).
\]
Hence, for any integer $m\ge 0$,
\[
\|D^m a_0\|_{L^2}
=\delta\lambda^{1-\beta}\|D^m(F(\lambda\cdot, \lambda\cdot))\|_{L^2}
=\delta\lambda^{1-\beta+m-1}\|D^mF\|_{L^2}
\lesssim \delta\lambda^{m-\beta}.
\]
Therefore, for any integer $k\ge 0$,
\[
\|a_0\|_{H^k}\lesssim \sum_{m=0}^k\|D^m a_0\|_{L^2}\lesssim \delta\lambda^{k-\beta}.
\]
For general $s\ge 0$, write $s=k+\sigma$ with $k\in\mathbb N_0$ and $\sigma\in[0,1)$.
Using interpolation between $H^k$ and $H^{k+1}$,
\[
\|a_0\|_{H^s}\lesssim \|a_0\|_{H^k}^{1-\sigma}\|a_0\|_{H^{k+1}}^{\sigma}
\lesssim \delta\lambda^{(k-\beta)(1-\sigma)+(k+1-\beta)\sigma}
=\delta\lambda^{s-\beta}.
\]
Similarly we have
\[
\begin{split}
\|u_0\|_{L^2}^2&=\delta^2\int_0^{2\pi}\int_0^\infty \lambda^{2(3-\beta)}(h'(\lambda r))^2r\, drd\theta\\
&\lesssim \delta^2\int_0^\infty  \lambda^{2(2-\beta)}(h'(\lambda r))^2(\lambda r)\, d(\lambda r)\\
&\lesssim \delta^2\lambda^{2(2-\beta)}
\end{split}
\]
and hence obtain the $H^s$ norm in \eqref{est-u0}. The claimed $C^1$ norm of $u_0$ is obvious as well.

\cbdu

\medskip

\subsection{Approximating solution}
\label{sec-app}

In polar coordinate, we can write
\[u_0\cdot\nabla \bar a=(u_0)_r\partial_r\bar a+\frac{(u_0)_\theta}{r}\partial_\theta \bar a=\frac{\partial_r b_0}{r}\partial_\theta \bar a\]
since $(u_0)_r=0$ and $(u_0)_\theta=\partial_r b_0$ in view of \eqref{initial-u}. Hence the transport equation \eqref{app-a} becomes
\[\partial_t \bar a+\frac{\partial_r b_0}{r}\partial_\theta \bar a=0 \]
which has the solution 
\begin{equation}\label{app-a-sol}
\bar a=\delta\lambda^{1-\beta} g(\lambda r)\cos\left(\theta- \frac{\partial_r b_0}{r}t\right).
\end{equation}
%Denote 
%\[
%\begin{split}
%\bar a_{rad}&:= \lambda^{1-\beta}f(\lambda r), \\ 
%\bar a_{osc}&:= \lambda^{1-\beta} g(\lambda r) \lambda^{-\beta\gamma}\cos\left(\lambda^\gamma(\theta- \frac{\partial_r b_0}{r}t)\right).
%\end{split}
%\]

\begin{Lemma}\label{le-bar-a-negative}
For $\eta\geq 0$, we have
\begin{equation}\notag
\|\bar a(t)\|_{\dot H^{-\eta}} \lesssim  \delta(1+\delta\lambda^{4-\beta}t)^{-\eta}\lambda^{-\eta-\beta}
\end{equation}
for any $t\in(0, T]$.
\end{Lemma}
\pf
First of all, we have 
\begin{equation}\label{osc-l2}
\begin{split}
\|\bar a\|_{L^2}^2&=\delta^2\int_0^{2\pi}\int_0^\infty \lambda^{2} g^2(\lambda r) \lambda^{-2\beta}\cos^2\left(\theta- \frac{\partial_r b_0}{r}t\right) r\, drd\theta\\
&\sim  \delta^2\int_0^\infty g^2(\lambda r) \lambda^{-2\beta} (\lambda r)\, d(\lambda r)\\
&\sim \delta^2\lambda^{-2\beta}.
\end{split}
\end{equation}
On the other hand, we note
\[\cos\left(\theta- \frac{\partial_r b_0}{r}t\right)=\cos\left(\theta-\delta\lambda^{4-\beta}\tilde h(\lambda r)t \right)\]
with $\tilde h(\lambda r)=\frac{h'(\lambda r)}{\lambda r}$. Note the derivative on $\tilde h$ gives a factor $\lambda$. 
Since $\supp g\subset(1,4)$ and $h'=1$ on $(2,3)$, on $\supp g(\lambda r)\cap (2,3)$ we have
\[
\bar a(t,r,\theta)=\delta\lambda^{1-\beta}g(\lambda r)\cos\!\left(\theta-\frac{\delta\lambda^{3-\beta}t}{r}\right).
\]
Define
\[
\omega(r,\theta):=\theta-\frac{\delta\lambda^{3-\beta}t}{r},
\qquad
\mu(t):=\Bigl(\lambda^2+(\delta\lambda^{5-\beta}t)^2\Bigr)^{1/2}.
\]
On $\supp g(\lambda r)\cap (2,3)$ (where $r\sim\lambda^{-1}$), one has
\[
\nabla\omega
=\frac{\delta\lambda^{3-\beta}t}{r^2}e_r+\frac1r e_\theta,
\qquad
|\nabla\omega|\sim \mu(t),
\]
and, for all integers $m\ge1$,
\[
|D^m\omega|\lesssim_m \mu(t)\lambda^{m-1}.
\]

Write $\bar a=\Re\!\big(Ae^{i\omega}\big)$ with $A(r):=\delta\lambda^{1-\beta}g(\lambda r)$.
Let $\psi\in C_c^\infty(\R^2)$ and define
\[
\mathcal L:=\frac{\nabla\omega}{i|\nabla\omega|^2}\cdot\nabla,
\qquad
\mathcal L(e^{i\omega})=e^{i\omega}.
\]
For any integer $m\ge0$, integrating by parts $m$ times gives
\[
\left|\int_{\R^2}Ae^{i\omega}\psi\,dx\right|
=\left|\int_{\R^2}e^{i\omega}(\mathcal L^*)^m(A\psi)\,dx\right|
\lesssim_m \mu(t)^{-m}\,\delta\lambda^{-\beta}\,\|\psi\|_{H^m}.
\]
Hence, by duality,
\[
\|\bar a(t)\|_{\dot H^{-m}}\lesssim_m \mu(t)^{-m}\,\delta\lambda^{-\beta}
\le (\delta\lambda^{5-\beta}t)^{-m}\delta\lambda^{-\beta},
\qquad m\in\mathbb N.
\]
Now fix $\eta\ge0$ and choose an integer $m>\eta$. Combining this with
\eqref{osc-l2} and interpolating between $\dot H^0$ and $\dot H^{-m}$ with
interpolation parameter $\eta/m\in(0,1)$ yields
\[
\|\bar a(t)\|_{\dot H^{-\eta}}
\lesssim \delta(\lambda+\delta\lambda^{5-\beta}t)^{-\eta}\lambda^{-\beta}\lesssim \delta(1+\delta\lambda^{4-\beta}t)^{-\eta}\lambda^{-\eta-\beta}.
\]

%Applying Lemma 8 from \cite{CMZO} we infer
%\begin{equation}\notag
%\|\bar a(t)\|_{\dot H^{-\eta}} \lesssim  (\delta\lambda^{5-\beta}t)^{-\eta}\|\bar a\|_{L^2}\lesssim \delta(\delta\lambda^{5-\beta}t)^{-\eta}\lambda^{-\beta}.
%\end{equation}
%{\color{blue} The derivative of $\cos( \theta)$ costs a factor $\lambda$. Assume $\lambda^{-\eta}\lesssim (\delta\lambda^{5-\beta}t)^{-\eta}$, which is equivalent to $\delta\lambda^{4-\beta}t \lesssim 1$. This is compatible with the choice of $t_N$ in Lemma \ref{le-norm-inflation}.}

\cbdu

\begin{Lemma}\label{le-a-osc}
For $s> 0$, we have for $t>0$
\begin{equation}\notag
\|\bar a(t)\|_{\dot H^{s}}\sim \delta(1+\delta\lambda^{4-\beta}t)^s\lambda^{s-\beta}.
%\ \ \|\bar a_{osc}(t)\|_{C^\zeta}\lesssim \lambda^{(5-\beta+\gamma)\zeta+1-\beta-\beta\gamma}.
\end{equation}
\end{Lemma}
\pf
Applying Lemma \ref{le-bar-a-negative}, \eqref{osc-l2} and interpolation yields
\begin{equation}\notag
\begin{split}
\|\bar a(t)\|_{\dot H^{s}} ^{\frac{\eta}{\eta+s}}&\gtrsim \|\bar a(t)\|_{L^2} \|\bar a(t)\|_{\dot H^{-\eta}} ^{-\frac{s}{\eta+s}}\\
&\gtrsim \delta\lambda^{-\beta}\left[\delta^{-1}(1+\delta\lambda^{4-\beta}t)^\eta \lambda^{\eta+\beta}\right]^{\frac{s}{\eta+s}}\\
&\sim \delta^{\frac{\eta}{\eta+s}}(1+\delta\lambda^{4-\beta}t)^{\frac{\eta s}{\eta+s}} \lambda^{(s-\beta)\frac{\eta}{\eta+s}}
\end{split}
\end{equation}
which gives the lower bound
\[\|\bar a(t)\|_{\dot H^{s}} \gtrsim  \delta(1+\delta\lambda^{4-\beta}t)^s\lambda^{s-\beta}.\]
On the other hand, the upper bound
\[
\|\bar a(t)\|_{\dot H^{s}} \lesssim  (\lambda+\delta\lambda^{5-\beta}t)^{s}\|\bar a(t)\|_{L^2}\lesssim  \delta(1+\delta\lambda^{4-\beta}t)^s\lambda^{s-\beta}
\]
is obvious. It completes the proof.
\cbdu

\begin{Lemma}\label{le-norm-inflation}
Let $t_N=\delta^{-\frac2\beta-1}\lambda^{\beta-4}$. We have for sufficiently small $\delta>0$
\begin{equation}\notag
\|\bar a(t_N)\|_{\dot H^{\beta}}\sim \frac1\delta.
\end{equation}
\end{Lemma}
\pf
It follows from Lemma \ref{le-a-osc} that 
\[ \|\bar a(t_N)\|_{\dot H^{\beta}}\sim \delta (1+\delta\lambda^{4-\beta}\delta^{-\frac2\beta-1}\lambda^{\beta-4})^\beta \lambda^{\beta-\beta}\sim \delta (1+\delta^{-\frac2\beta})^\beta \sim \delta^{-1}.\]

\cbdu

%Note for $\beta<4$, $0<\zeta<4-\beta$ and $s\geq \beta$, we have $(5-\beta+\gamma-\zeta)s-\beta\gamma-\beta>s-\beta$ and hence
%\begin{equation}\label{est-a-bar}
%\|\bar a(t)\|_{\dot H^{s}}\lesssim \lambda^{(5-\beta+\gamma-\zeta)s-\beta\gamma}, \ \ \forall \ t\in[0,t_N].
%\end{equation}

\begin{Remark}\label{rk-bar-a}
From now on we apply Lemma \ref{le-a-osc} as 
\[\|\bar a(t)\|_{\dot H^{s}}\sim \delta(\delta\lambda^{4-\beta}t)^s\lambda^{s-\beta}\]
with the understanding that we work with the time scale $t\sim \delta\lambda^{4-\beta}$ away from 0 and hence 
\[1+\delta\lambda^{4-\beta}t\sim \delta\lambda^{4-\beta}t.\]
\end{Remark}

\begin{Lemma}\label{le-bar-u}
We have for $0< t\leq t_N$ and $0<s\leq \beta$, 
\begin{equation}\notag
\|\bar u(t)\|_{H^{s-2}}\lesssim \delta \lambda^{s-\beta}+\delta^{1-\frac{2s+10}{\beta}}\lambda^{s-\beta}.
\end{equation}
%In particular, 
%\[\|\bar u(t_N)\|_{H^{\beta-2}}\leq \|u_0\|_{H^{\beta-2}}+C\delta^{-\frac{10}\beta}.\]
\end{Lemma}
\pf
According to \eqref{app-u}, we have
\[
\bar u(t)=u_0+\int_0^t \nabla^{\perp}(\nabla^{\perp}\bar a(\tau)\cdot\nabla\Delta\bar a(\tau)) \, d\tau.
\]
Thus
\begin{equation}\notag
\begin{split}
\|\bar u(t)\|_{L^2}\leq& \ \|u_0\|_{L^2}+t\|\nabla^{\perp}(\nabla^{\perp}\bar a\cdot\nabla\Delta\bar a)\|_{L^2}
\end{split}
\end{equation}
%Recall \[\|\bar a_{rad}\|_{H^\eta}\lesssim \lambda^{\eta-\beta}, \ \ \ \|\bar a_{rad}\|_{C^\eta}\lesssim \lambda^{1+\eta-\beta}.\]
Applying H\"older's inequality and Gagliardo-Nirenberg's inequality gives
\begin{equation}\notag
\begin{split}
&\|\nabla^{\perp}(\nabla^{\perp}\bar a\cdot\nabla\Delta\bar a)\|_{L^2}\\
\leq&\ \|D^2\bar a\|_{L^\infty} \|\nabla\Delta \bar a\|_{L^2}+\|D\bar a\|_{L^\infty} \|D^4 \bar a\|_{L^2}\\
\leq&\ \|D^2\bar a\|_{H^1} \|\nabla\Delta \bar a\|_{L^2}+\|D\bar a\|_{H^1} \|D^4 \bar a\|_{L^2}.
%\lesssim&\ \delta(\delta\lambda^{4-\beta}t)^3\lambda^{3-\beta}\delta(\delta\lambda^{4-\beta}t)^3\lambda^{3-\beta}+%\delta(\delta\lambda^{4-\beta}t)^2\lambda^{2-\beta}\delta(\delta\lambda^{4-\beta}t)^4\lambda^{4-\beta}\\
%\lesssim&\ \delta^8\lambda^{30-8\beta}t^6.
\end{split}
\end{equation}
Therefore, we infer, using Lemma \ref{le-a-osc} and Remark \ref{rk-bar-a}
\begin{equation}\notag
\begin{split}
&\|\nabla^{\perp}(\nabla^{\perp}\bar a\cdot\nabla\Delta\bar a)\|_{L^2}\\
\lesssim&\ \delta(\delta\lambda^{4-\beta}t)^3\lambda^{3-\beta}\delta(\delta\lambda^{4-\beta}t)^3\lambda^{3-\beta}+
\delta(\delta\lambda^{4-\beta}t)^2\lambda^{2-\beta}\delta(\delta\lambda^{4-\beta}t)^4\lambda^{4-\beta}\\
\lesssim&\ \delta^8\lambda^{30-8\beta}t^6.
\end{split}
\end{equation}

Summarizing the estimates above yields 
\begin{equation}\label{bar-u-tN}
\begin{split}
\|\bar u(t)\|_{L^2}\lesssim& \ \|u_0\|_{L^2}+\delta^8\lambda^{30-8\beta}t^7, \quad  0< t\leq t_N.
\end{split}
\end{equation}
We note a derivative on $\bar a(t)$ costs a factor of $\delta\lambda^{5-\beta}t$. Thus we infer from \eqref{bar-u-tN} for $0< t\leq t_N$,
\begin{equation}\label{bar-u-tN-2}
\begin{split}
\|\bar u(t)\|_{H^{s-2}}\leq& \ \|u_0\|_{H^{s-2}}
+C\delta^8(\delta\lambda^{5-\beta}t)^{s-2}\lambda^{30-8\beta}t^7
%&+C(\lambda^{5-\beta+\gamma}t_N)^{\beta-2}\lambda^{6(5-\beta+\gamma)-2\beta-2\beta\gamma}t_N^7
\end{split}
\end{equation}
for some constant $C>0$.  Recall $t_N=\delta^{-\frac2\beta-1}\lambda^{\beta-4}$. It follows from \eqref{bar-u-tN-2} that
\begin{equation}\label{bar-u-tN-s}
\begin{split}
\|\bar u(t)\|_{H^{s-2}}\leq& \ C\delta\lambda^{s-\beta}+C\delta^{1-\frac{2s+10}{\beta}}\lambda^{s-\beta}, \quad 0< t\leq t_N
\end{split}
\end{equation}
and in particular
\begin{equation}\label{bar-u-tN-3}
\begin{split}
\|\bar u(t_N)\|_{H^{\beta-2}}\lesssim& \ \delta+\delta^{-1-\frac{10}\beta}.
\end{split}
\end{equation}
%Take $T=C^{-1}\lambda^{-\gamma(s-2)-4\gamma+\beta\gamma+6\beta-22}$. When $s=\beta$, $T=C^{-1}\lambda^{-2\gamma+6\beta-22}$.
\cbdu

\bigskip

\section{Control of perturbation}
\label{sec-perturb}

Let $(a, u)$ be a solution of \eqref{emhd-au} with the initial data \eqref{data}-\eqref{initial-u}.
Denote the perturbation $A=a-\bar a$ which satisfies 
\begin{equation}\label{eq-A}
\begin{split}
\partial_tA+u\cdot\nabla A+(u-u_0)\cdot\nabla\bar a=&\ 0,\\
A(x,0)=&\ 0.
\end{split}
\end{equation}
On the other hand, it follows from the second equation of \eqref{emhd-au} that
\begin{equation}\label{u-int}
\begin{split}
u(t)-u_0=&\int_0^t\nabla^{\perp}(\nabla^{\perp}a(\tau)\cdot\nabla\Delta a(\tau))\,d\tau.
%=&\int_0^t\nabla^{\perp}(\nabla^{\perp}A(\tau)\cdot\nabla\Delta A(\tau))\,d\tau+\int_0^t\nabla^{\perp}(\nabla^{\perp}A(\tau)\cdot\nabla\Delta \bar a(\tau))\,d\tau\\
%&+\int_0^t\nabla^{\perp}(\nabla^{\perp}\bar a(\tau)\cdot\nabla\Delta A(\tau))\,d\tau+\int_0^t\nabla^{\perp}(\nabla^{\perp}\bar a(\tau)\cdot\nabla\Delta \bar a(\tau))\,d\tau.
\end{split}
\end{equation}
We observe that for $t>0$ small enough, % (depending on the norms of $\bar a$ and $A$), 
$u(t)-u_0$ is expected to be controlled and so is $A(t)$ according to \eqref{eq-A}.  In particular, the norm $\|u(t)-u_0\|_{L^2}$ depends on the norm $\|\nabla a\|_{H^3}$, and higher norm of $u(t)-u_0$ depends on higher norm $\|\nabla a\|_{H^s}$ with $s\geq 3$. For this purpose, we first establish the following higher norm estimate of $a$ by applying a continuity argument.

\medskip

\subsection{Higher norm estimate of $a$}
\label{sec-high-norm}

\begin{Lemma}\label{le-high}
Let $s\geq 4$. There exists a constant $M>0$ such that 
\begin{equation}\label{high-norm}
\|a(t)\|_{H^s}\leq 2M\delta\lambda^{s-\beta}
\end{equation}
for all $t\in[0, t_N]$.
\end{Lemma}
\pf
Due to \eqref{est-u0}, we have at the initial time 
\[\|a_0\|_{H^s}\lesssim \delta\lambda^{s-\beta}\leq C \delta\lambda^{s-\beta}\]
for a constant $C>0$. 
Thus there exist a (large) constant $M>0$ and a small enough time $t_0<t_N$ such that
\begin{equation}\label{est-a-hs-t0}
\|a(t)\|_{H^s}\leq 2M\delta\lambda^{s-\beta}, \ \ t\in[0,t_0]
\end{equation}
for $s\geq 4$. Then we can show that the estimate holds on $[0,t_N]$ using a bootstrapping argument. For this purpose, we 
fix $s\geq 4$, and define the maximum time of the bootstrap bound as
\[
T_*:=\sup\Bigl\{t\in[0,t_N]: \|a(\tau)\|_{H^s}\le 2M\delta\lambda^{s-\beta}\ \ \forall\,\tau\in[0,t]\Bigr\}.
\]
It suffices to prove $T_*=t_N$.

For any $t\in(0,T_*]$, the bootstrap bound implies
\[
\|a(\tau)\|_{H^r}\le 2M\delta\lambda^{r-\beta},\qquad 4\le r\le s,\ \ 0\le\tau\le t.
\]
Also, from the transport feature of the $a$ equation in \eqref{emhd-au} and the basic energy law we have for $0\leq t\leq t_N$
\begin{equation}\label{est-basic-a}
\|a(\tau)\|_{L^2}\le \|a_0\|_{L^2}\lesssim \delta\lambda^{-\beta},\qquad \|Da(\tau)\|_{L^2}\le \|Da_0\|_{L^2}+\|b_0\|_{L^2}\lesssim \delta\lambda^{1-\beta}.
\end{equation}

We deduce from \eqref{u-int} using H\"older's inequality and Gagliardo-Nirenberg's inequality, for $\alpha\geq 0$
\begin{equation}\label{u-u0-basic}
\begin{split}
\|u(t)-u_0\|_{H^\alpha}\leq & \int_0^t\|\nabla^{\perp}(\nabla^{\perp}a(\tau)\cdot\nabla\Delta a(\tau))\|_{H^\alpha}\, d\tau\\
\leq & \int_0^t\sum_{\alpha'=0}^{\alpha+1}\|D^{\alpha'}\nabla^{\perp}a(\tau)D^{\alpha+1-\alpha'}\nabla\Delta a(\tau)\|_{L^2}\, d\tau\\
\leq & \int_0^t\sum_{\alpha'=0}^{\alpha+1}\|D^{\alpha'}\nabla^{\perp}a(\tau)\|_{L^\infty}\|D^{\alpha+1-\alpha'}\nabla\Delta a(\tau)\|_{L^2}\, d\tau\\
\leq & \int_0^t\sum_{\alpha'=0}^{\alpha+1}\|a(\tau)\|_{H^{2+\alpha'}}\|a(\tau)\|_{H^{4+\alpha-\alpha'}}\, d\tau\\
\lesssim& \int_0^t\|Da(\tau)\|_{L^2}^{\frac{2+\alpha}{3+\alpha}}\|D a(\tau)\|_{H^{3+\alpha}}^{\frac{4+\alpha}{3+\alpha}} d\tau.
\end{split}
\end{equation}
%Applying the basic energy estimate we have 
%\[\|Da(t)\|_{L^2}\leq \|Da_0\|_{L^2}+\|b_0\|_{L^2}\lesssim \delta\lambda^{1-\beta}.\]
Hence combined with the estimate \eqref{est-a-hs-t0}, it follows for $0<t\leq t_0$
\begin{equation}\label{u-H}
\|u(t)-u_0\|_{H^\alpha}\leq C_0(2M\delta)^{\frac{4+\alpha}{3+\alpha}}t_0\delta^{\frac{2+\alpha}{3+\alpha}}\lambda^{\frac{2+\alpha}{3+\alpha}(1-\beta)}\lambda^{\frac{4+\alpha}{3+\alpha}(4+\alpha-\beta)}
\end{equation}
for a constant $C_0>0$. We note that the gaining factor from the estimate of $\|u(t)-u_0\|_{H^\alpha}$ to $\|u(t)-u_0\|_{H^{\alpha+1}}$ is 
\begin{equation}\label{u-H1}
\begin{split}
&\lambda^{(\frac{3+\alpha}{4+\alpha}-\frac{2+\alpha}{3+\alpha})(1-\beta)}\lambda^{\frac{5+\alpha}{4+\alpha}((5+\alpha)-\beta)-\frac{4+\alpha}{3+\alpha}((4+\alpha)-\beta)}=\lambda.
\end{split}
\end{equation}

It is obvious that for $\alpha=0,1$ we have the estimates in \eqref{est-basic-a}.
We further deduce from \eqref{emhd-au} for $\alpha\geq 2$
%\begin{equation}\notag
%\partial_t D^\alpha a+D^\alpha(u\cdot\nabla a)=0
%\end{equation}
\begin{equation}\notag
\frac12\frac{d}{dt}\|D^\alpha a\|_{L^2}^2+\int_{\mathbb R^2} D^\alpha(u\cdot\nabla a) D^\alpha a\, dxdy=0.
\end{equation}
It follows
\begin{equation}\notag
\frac12\frac{d}{dt}\|D^\alpha a\|_{L^2}^2\leq \sum_{\alpha'=0}^{\alpha-1}\|D^{\alpha-\alpha'}u\|_{L^\infty}\|D^{\alpha'}\nabla a\|_{L^2}\|D^{\alpha} a\|_{L^2}
\end{equation}
and thus
\begin{equation}\notag
\begin{split}
\frac{d}{dt}\|D^\alpha a\|_{L^2}\leq& \sum_{\alpha'=0}^{\alpha-1}\|D^{\alpha-\alpha'}u\|_{L^\infty}\|D^{\alpha'}\nabla a\|_{L^2}\\
\leq& \sum_{\alpha'=0}^{\alpha-1}\|\nabla u\|_{H^{\alpha-\alpha'}}\|\nabla a\|_{H^{\alpha'}}\\
=& \sum_{\alpha'=0}^{\alpha-2}\|\nabla u\|_{H^{\alpha-\alpha'}}\|\nabla a\|_{H^{\alpha'}}+\|\nabla u\|_{H^1}\|\nabla a\|_{H^{\alpha}}.
\end{split}
\end{equation}
Applying Gr\"onwall's inequality we obtain
\begin{equation}\label{est-a-H}
\begin{split}
\|D^\alpha a(t)\|_{L^2}\leq&\ \|D^\alpha a_0\|_{L^2}e^{\int_0^t\|\nabla u(\tau)\|_{H^1}\, d\tau}\\
&+\int_0^t e^{\int_{\tau}^t\|\nabla u(\tau')\|_{H^1}\, d\tau'}  \sum_{\alpha'=0}^{\alpha-2}\|\nabla u(\tau)\|_{H^{\alpha-\alpha'}}\|\nabla a(\tau)\|_{H^{\alpha'}}\, d\tau.
\end{split}
\end{equation}
Taking $\alpha=2$ in \eqref{u-H} gives
\begin{equation}\notag
\|u(t)-u_0\|_{H^2}\leq C_0\delta^2(2M)^{\frac65}t_0\lambda^{\frac45(1-\beta)}\lambda^{\frac65(6-\beta)}.
\end{equation}
In view of $\|u_0\|_{H^2}\lesssim \delta\lambda^{4-\beta}$ from \eqref{est-u0},
we can choose 
\begin{equation}\label{t0-choice1}
t_0<\frac12C_0^{-1}(2M)^{-\frac65}\delta^{-1}\lambda^{\beta-4}, \ \ t_0<\delta^{-1}\lambda^{\beta-4} C^{-1}\ln 2
\end{equation}
 such that for $0\leq t\leq \min\{t_0, T_*\}$
\begin{equation}\notag
\|u(t)-u_0\|_{H^2}\leq C_0\delta^2(2M)^{\frac65}t_0\lambda^{\frac45(1-\beta)}\lambda^{\frac65(6-\beta)}< \frac12\delta \lambda^{4-\beta}
\end{equation}
and
\begin{equation}\notag
\begin{split}
e^{\int_0^t\|\nabla u(\tau)\|_{H^1}\, d\tau}\leq e^{t_0(\|u(t)-u_0\|_{H^2}+\|u_0\|_{H^2})}\leq e^{Ct_0\delta\lambda^{4-\beta}}<2.
\end{split}
\end{equation}
It then follows from \eqref{est-a-H} that
\begin{equation}\label{est-a-H-general}
\begin{split}
\|D^\alpha a(t)\|_{L^2}\leq&\ 2\|D^\alpha a_0\|_{L^2}+2\int_0^t \sum_{\alpha'=2}^{\alpha}\|\nabla u(\tau)\|_{H^{\alpha'}}\|\nabla a(\tau)\|_{H^{\alpha-\alpha'}}\, d\tau\\
\leq&\ 2\|D^\alpha a_0\|_{L^2}+2t_0 \sum_{\alpha'=2}^{\alpha}\|\nabla u(t)\|_{H^{\alpha'}}\|\nabla a(t)\|_{H^{\alpha-\alpha'}}\\
%&+2t_0 \sum_{\alpha'=0}^{\alpha-2}\|\nabla u_0\|_{H^{\alpha-\alpha'}}\|\nabla a(t)\|_{H^{\alpha'}}
\end{split}
\end{equation}
which is an iterative estimate of $\|D^\alpha a(t)\|_{L^2}$ for $\alpha\geq 2$ on $[0, \min\{t_0,T_*\}]$. Iterating \eqref{est-a-H-general} from $\alpha=2$ to $\alpha=4$ yields
\begin{equation}\label{est-a-H2}
\begin{split}
\|D^2 a(t)\|_{L^2}\leq&\ 2\|D^2 a_0\|_{L^2}+2t \|\nabla u(t)\|_{H^2}\|\nabla a(t)\|_{L^2},\\
\|D^3 a(t)\|_{L^2}\leq&\ 2\|D^3 a_0\|_{L^2}+2t \|\nabla u(t)\|_{H^3}\|\nabla a(t)\|_{L^2}\\
&+2t \|\nabla u(t)\|_{H^2}\|\nabla a(t)\|_{H^1},\\
\|D^4 a(t)\|_{L^2}\leq&\ 2\|D^4 a_0\|_{L^2}+2t \|\nabla u(t)\|_{H^4}\|\nabla a(t)\|_{L^2}\\
&+2t \|\nabla u(t)\|_{H^3}\|\nabla a(t)\|_{H^1}+2t \|\nabla u(t)\|_{H^2}\|\nabla a(t)\|_{H^2}.
\end{split}
\end{equation}
Since $\|\nabla a(t)\|_{L^2}\leq\|\nabla a_0\|_{L^2}+\|b_0\|_{L^2} \lesssim \delta\lambda^{1-\beta}$, the structure of inequalities in \eqref{est-a-H2} with $t$ factor in front of the quadratic terms indicates that there exists $t_0\ll 1$ depending only on $\beta$ such that for $t\in [0, \min\{t_0,T_*\}]$
\begin{equation}\label{est-a-rough}
\|D^2 a(t)\|_{L^2}\leq 4\|D^2 a_0\|_{L^2}, \ \ \|D^3 a(t)\|_{L^2}\leq 4\|D^3 a_0\|_{L^2}
\end{equation}
and 
\begin{equation}\label{est-a-H4}
\|D^4a(t)\|_{L^2}\leq M\delta\lambda^{4-\beta}.
\end{equation}
Now we assume for $4\leq s\leq \alpha-1$ the estimate
\begin{equation}\label{high-norm-better}
\|D^sa(t)\|_{L^2}\leq M\delta\lambda^{s-\beta},  \qquad t\in [0, \min\{t_0,T_*\}]
\end{equation}
is satisfied. The goal is to show \eqref{high-norm-better} holds for $s=\alpha$. Rewriting \eqref{est-a-H-general} gives
\begin{equation}\label{est-a-H5}
\begin{split}
\|D^\alpha a(t)\|_{L^2}\leq&\ 2\|D^\alpha a_0\|_{L^2}+2t \sum_{\alpha'=2}^{\alpha}\|\nabla (u(t)-u_0)\|_{H^{\alpha'}}\|\nabla a(t)\|_{H^{\alpha-\alpha'}}\\
&+2t \sum_{\alpha'=2}^{\alpha}\|\nabla u_0\|_{H^{\alpha'}}\|\nabla a(t)\|_{H^{\alpha-\alpha'}}.
\end{split}
\end{equation}
%Note when $\alpha=2$, it is 
%\begin{equation}\notag
%\|D^2 a(t)\|_{L^2}\leq 2\|D^2 a_0\|_{L^2}+2t_0 \|\nabla u(t)\|_{H^2}\|\nabla a(t)\|_{L^{2}};
%\end{equation}
%while for $\alpha\geq 3$ we have
%\begin{equation}\notag
%\begin{split}
%\|D^\alpha a(t)\|_{L^2}
%\leq&\ 2\|D^\alpha a_0\|_{L^2}+2t_0\|\nabla u(t)\|_{H^2}\|\nabla a(t)\|_{H^{\alpha-2}}\\
%&+2t_0 \sum_{\alpha'=0}^{\alpha-3}\|\nabla u(t)\|_{H^{\alpha-\alpha'}}\|\nabla a(t)\|_{H^{\alpha'}}.
%\end{split}
%\end{equation}
%We observe that the gaining factor of estimate from $\|D^\alpha a\|_{L^2}$ to $\|D^{\alpha+1} a\|_{L^2}$ is $t_0\|\nabla u(t)\|_{H^2}$. 
First of all, it follows from \eqref{est-u0} 
\begin{equation}\label{est-middle-1}
2\|D^\alpha a_0\|_{L^2}\leq C \delta\lambda^{\alpha-\beta}<\frac{M}{10}\delta\lambda^{\alpha-\beta}
\end{equation}
 for some large $M>0$.

Combining \eqref{u-H} and \eqref{high-norm-better} gives
\begin{equation}\notag
\begin{split}
&2t \sum_{\alpha'=2}^{\alpha}\|\nabla (u(t)-u_0)\|_{H^{\alpha'}}\|\nabla a(t)\|_{H^{\alpha-\alpha'}}\\
\leq&\ 2C_0Mt^2 \sum_{\alpha'=2}^{\alpha-3}(2M)^{\frac{5+\alpha'}{4+\alpha'}}\delta^2\lambda^{\frac{3+\alpha'}{4+\alpha'}(1-\beta)}\lambda^{\frac{5+\alpha'}{4+\alpha'}((5+\alpha')-\beta)}
\delta\lambda^{(\alpha+1-\alpha')-\beta}\\
&+2t\|\nabla (u(t)-u_0)\|_{H^{\alpha-2}}\|\nabla a(t)\|_{H^{2}}\\
&+2t\|\nabla (u(t)-u_0)\|_{H^{\alpha-1}}\|\nabla a(t)\|_{H^{1}}+2t\|\nabla (u(t)-u_0)\|_{H^{\alpha}}\|\nabla a(t)\|_{L^{2}}.
\end{split}
\end{equation}
We deduce
\begin{equation}\notag
\begin{split}
& 2C_0Mt^2 \sum_{\alpha'=2}^{\alpha-3}(2M)^{\frac{5+\alpha'}{4+\alpha'}}\delta^2\lambda^{\frac{3+\alpha'}{4+\alpha'}(1-\beta)}\lambda^{\frac{5+\alpha'}{4+\alpha'}((5+\alpha')-\beta)}\\
&\cdot\delta\lambda^{(\alpha+1-\alpha')-\beta}\\
\leq&\ 4C_0Mt^2(2M)^{\frac{7}{6}}\delta^3\lambda^{\frac{5}{6}(1-\beta)}\lambda^{\frac{7}{6}(7-\beta)}\lambda^{(\alpha-1)-\beta}\\
\leq&\ 4C_0t^2(2M)^{3}\delta^3\lambda^{2(4-\beta)}\lambda^{\alpha-\beta},
\end{split}
\end{equation}
thanks to \eqref{u-H1} and the fact that the derivative cost of $u(t)-u_0$ and $a(t)$ is both order $\lambda$. Using \eqref{est-u0}, \eqref{u-H} and \eqref{est-a-rough}, we obtain
\begin{equation}\notag
\begin{split}
&2t\|\nabla (u(t)-u_0)\|_{H^{\alpha-2}}\|\nabla a(t)\|_{H^{2}}+2t\|\nabla (u(t)-u_0)\|_{H^{\alpha-1}}\|\nabla a(t)\|_{H^{1}}\\
&+2t\|\nabla (u(t)-u_0)\|_{H^{\alpha}}\|\nabla a(t)\|_{L^{2}}\\
\leq&\ 2 C_0t^2  (2M)^{\frac{3+\alpha}{2+\alpha}}\delta^2\lambda^{\frac{1+\alpha}{2+\alpha}(1-\beta)}\lambda^{\frac{3+\alpha}{2+\alpha}((3+\alpha)-\beta)}\delta\lambda^{3-\beta}\\
&+2C_0t^2 (2M)^{\frac{4+\alpha}{3+\alpha}}\delta^2\lambda^{\frac{2+\alpha}{3+\alpha}(1-\beta)}\lambda^{\frac{4+\alpha}{3+\alpha}((4+\alpha)-\beta)}\delta\lambda^{2-\beta}\\
&+2C_0t^2(2M)^{\frac{5+\alpha}{4+\alpha}}\delta^2\lambda^{\frac{3+\alpha}{4+\alpha}(1-\beta)}\lambda^{\frac{5+\alpha}{4+\alpha}((5+\alpha)-\beta)}\delta\lambda^{1-\beta}\\
\leq&\ 8C_0t^2(2M)^{\frac{3+\alpha}{2+\alpha}}\delta^3\lambda^{\frac{3+\alpha}{4+\alpha}(1-\beta)}\lambda^{\frac{5+\alpha}{4+\alpha}((5+\alpha)-\beta)}\lambda^{1-\beta}\\
\leq&\ 8C_0t^2(2M)^{2}\delta^3\lambda^{2(4-\beta)} \lambda^{\alpha-\beta}.
\end{split}
\end{equation}
We choose small $t_0$ satisfying
\begin{equation}\label{choice-t-2}
\begin{split}
0<t_0&<\frac{1}{8M\sqrt{10C_0}}\delta^{-1}\lambda^{\beta-4}.
%0<t_0&<\frac{1}{8\sqrt{5C_0M}}\lambda^{\beta-4}.
\end{split}
\end{equation}
Under the condition \eqref{choice-t-2}, we can verify 
\begin{equation}\notag
\begin{split}
4C_0t_0^2(2M)^{3}\delta^3\lambda^{2(4-\beta)}\lambda^{\alpha-\beta}
< \frac{M}{10}\delta\lambda^{\alpha-\beta}
\end{split}
\end{equation}
and
\begin{equation}\notag
\begin{split}
8C_0t_0^2(2M)^{2}\delta^3\lambda^{2(4-\beta)} \lambda^{\alpha-\beta}
< \frac{M}{10}\delta\lambda^{\alpha-\beta}.
\end{split}
\end{equation}
Summarizing the analysis above we have obtained
\begin{equation}\label{est-middle-2}
2t \sum_{\alpha'=2}^{\alpha}\|\nabla (u(t)-u_0)\|_{H^{\alpha'}}\|\nabla a(t)\|_{H^{\alpha-\alpha'}}\leq \frac{M}{5}\delta\lambda^{\alpha-\beta},  \quad t\in [0, \min\{t_0,T_*\}].
\end{equation}

Next we estimate the last term in \eqref{est-a-H5}, using \eqref{est-u0} and \eqref{high-norm-better}
\begin{equation}\notag
\begin{split}
&2t \sum_{\alpha'=2}^{\alpha}\|\nabla u_0\|_{H^{\alpha'}}\|\nabla a(t)\|_{H^{\alpha-\alpha'}}\\
\leq&\ 2Mt \sum_{\alpha'=2}^{\alpha}\delta\lambda^{\alpha'+3-\beta}\delta\lambda^{(\alpha+1-\alpha')-\beta}\\
\leq&\ 4Mt\delta^2\lambda^{5-\beta}\lambda^{(\alpha-1)-\beta}\\
=&\ 4Mt\delta^2\lambda^{4-\beta}\lambda^{\alpha-\beta}.
\end{split}
\end{equation} 
Thus by requiring
\begin{equation}\notag
t_0<\frac{1}{40}\delta^{-1}\lambda^{\beta-4}<t_N
\end{equation} 
which is satisfied in view of \eqref{choice-t-2} for large $M$,
it is clear to see 
$4Mt_0\delta^2\lambda^{4-\beta}<\frac{M}{10}\delta$, and hence
\begin{equation}\label{est-middle-3}
2t \sum_{\alpha'=2}^{\alpha}\|\nabla u_0\|_{H^{\alpha'}}\|\nabla a(t)\|_{H^{\alpha-\alpha'}}<
\frac{M}{10}\delta\lambda^{\alpha-\beta},  \qquad t\in [0, \min\{t_0,T_*\}].
\end{equation}

Therefore, we conclude that for $t_0$ small enough such that \eqref{est-a-rough} and \eqref{est-a-H4} are satisfied and \eqref{choice-t-2} is satisfied, combining \eqref{est-a-H5}, \eqref{est-middle-1}, \eqref{est-middle-2} and \eqref{est-middle-3}, the estimate \eqref{high-norm-better} holds for $s=\alpha$.
By induction, \eqref{high-norm-better} holds for all $s\geq 4$. In particular, the bootstrap constant $2M$ in \eqref{est-a-hs-t0} is improved on the time interval $[0, \min\{t_0,T_*\}]$.

To conclude globally on $[0,t_N]$, define
\[
t_j:=jt_0,\qquad j=0,1,\dots,J,\qquad J:=\left\lceil \frac{t_N}{t_0}\right\rceil .
\]
We claim
\[
\|a(t)\|_{H^s}\le 2M\delta\lambda^{s-\beta},
\qquad 0\le t\le \min\{t_j,T_*\},
\]
for every $j=0,\dots,J$. The case $j=1$ is exactly the estimate above on
$[0,\min\{t_0,T_*\}]$.

Assume the claim holds for some $j<J$. Repeat the same argument on the shifted
interval $[t_j,t_j+t_0]$, with initial time $t_j$ and initial data $a(t_j)$.
The bounds used in \eqref{est-middle-1}--\eqref{est-middle-3} depend only on
$\beta,\delta,\lambda,M$ and the bootstrap size $2M$, so the same $t_0$ works.
Hence
\[
\|a(t)\|_{H^s}\le M\delta\lambda^{s-\beta},
\qquad t_j\le t\le \min\{t_j+t_0,T_*\},
\]
which in particular implies the bootstrap bound $2M$ up to
$\min\{t_{j+1},T_*\}$. This proves the claim for all $j$ by induction.

Therefore $\|a(t)\|_{H^s}\le 2M\delta\lambda^{s-\beta}$ on $[0,T_*]$.
If $T_*<t_N$, pick $j$ with $T_*\in[t_j,t_{j+1}]$. The previous estimate gives
\[
\|a(T_*)\|_{H^s}\le M\delta\lambda^{s-\beta}<2M\delta\lambda^{s-\beta},
\]
and continuity of $t\mapsto\|a(t)\|_{H^s}$ extends the bootstrap bound slightly
beyond $T_*$, contradicting the definition of $T_*$. Thus $T_*=t_N$, and
\eqref{high-norm} follows on $[0,t_N]$.

\cbdu

\medskip

\subsection{Estimates of perturbation}
\label{sec-error}

We proceed to show the perturbations are under control.
\begin{Lemma}\label{le-perturbation}
Let $1<\beta<4$. We have for $s\leq \beta$
\[\|A(t)\|_{H^s}\lesssim \delta \lambda^{s-\beta}, \ \ \|u(t)-u_0\|_{H^{s-2}}\lesssim \delta^{1-\frac{2s+10}\beta}\lambda^{s-\beta}, \ \ \forall \ t\in[0,t_N]\]
\end{Lemma}
\pf
We first show the estimate on the short time interval $[0, t_0]$. 
 Multiplying the first equation of \eqref{eq-A} by $A$ and integrating over $\mathbb R^2$ yields
\begin{equation}\notag
\frac12\frac{d}{dt}\|A(t)\|_{L^2}^2\leq \|u(t)-u_0\|_{L^2}\|\nabla \bar a(t)\|_{L^\infty}\|A(t)\|_{L^2}
\end{equation}
where we used the fact $\nabla\cdot u=0$. It follows immediately
\begin{equation}\label{A-l2}
\|A(t)\|_{L^2}\leq 2\int_0^t\|u(\tau)-u_0\|_{L^2}\|\nabla \bar a(\tau)\|_{H^1}\, d\tau.
\end{equation}
%Thanks to Lemma \ref{le-high}, the estimate \eqref{u-H} holds on $[0,t_N]$. 
Applying \eqref{u-H} with $\alpha=0$ and Lemma \ref{le-a-osc} to \eqref{A-l2}, we infer for $0< t\leq t_0$
\begin{equation}\label{A-l2-final}
\begin{split}
\|A(t)\|_{L^2}&\leq 2C_0(2M)^{\frac43}\delta^2t_0^2\lambda^{\frac23(1-\beta)+\frac43\left(4-\beta\right)}\delta(\delta\lambda^{4-\beta}t_0)^2\lambda^{2-\beta}\\
&\leq CM^{-2}\delta \lambda^{-\beta}
\end{split}
\end{equation}
where we also used \eqref{choice-t-2}.

For $\alpha\geq 1$, acting $D^\alpha$ on \eqref{eq-A}, multiplying by $D^\alpha A$ and integrating over $\mathbb R^2$ we obtain
\begin{equation}\notag
\begin{split}
\frac12\frac{d}{dt}\|D^\alpha A\|_{L^2}^2&+\int_{\mathbb R^2} D^\alpha(u\cdot\nabla A)D^\alpha A\, dxdy\\
&+\int_{\mathbb R^2} D^\alpha((u-u_0)\cdot\nabla \bar a)D^\alpha A\, dxdy=0.
\end{split}
\end{equation}
Using H\"older's inequality and Gagliardo-Nirenberg yields
\begin{equation}\notag
\begin{split}
&\left|\int_{\mathbb R^2} D^\alpha(u\cdot\nabla A)D^\alpha A\, dxdy \right|\\
\leq &\sum_{\alpha'=1}^\alpha \left| \int_{\mathbb R^2}D^{\alpha'}u\cdot\nabla D^{\alpha-\alpha'}A D^\alpha A\, dxdy \right|\\
\leq &\sum_{\alpha'=1}^\alpha \|D^{\alpha'}u\|_{L^\infty}\|\nabla D^{\alpha-\alpha'}A\|_{L^2}\|D^\alpha A\|_{L^2}\\
\lesssim &\sum_{\alpha'=1}^\alpha \|\nabla u\|_{H^{\alpha'}}\|\nabla A\|_{H^{\alpha-\alpha'}}\|D^\alpha A\|_{L^2},
\end{split}
\end{equation}
and 
\begin{equation}\notag
\begin{split}
&\left|\int_{\mathbb R^2} D^\alpha((u-u_0)\cdot\nabla \bar a)D^\alpha A\, dxdy \right|\\
\leq &\sum_{\alpha'=0}^\alpha \left| \int_{\mathbb R^2}D^{\alpha'}(u-u_0)\cdot\nabla D^{\alpha-\alpha'}\bar a D^\alpha A\, dxdy \right|\\
\leq &\sum_{\alpha'=0}^\alpha \|D^{\alpha'}(u-u_0)\|_{L^2}\|\nabla D^{\alpha-\alpha'}\bar a\|_{L^\infty}\|D^\alpha A\|_{L^2}\\
\lesssim &\sum_{\alpha'=0}^\alpha \|u-u_0\|_{H^{\alpha'}}\|\nabla \bar a\|_{H^{\alpha-\alpha'+1}}\|D^\alpha A\|_{L^2}.
\end{split}
\end{equation}
Therefore we claim
\begin{equation}\notag
\begin{split}
\frac{d}{dt}\|D^\alpha A\|_{L^2}\lesssim \sum_{\alpha'=1}^\alpha \|\nabla u\|_{H^{\alpha'}}\|\nabla A\|_{H^{\alpha-\alpha'}}
+\sum_{\alpha'=0}^\alpha \|u-u_0\|_{H^{\alpha'}}\|\nabla \bar a\|_{H^{\alpha-\alpha'+1}}.
\end{split}
\end{equation}
Applying Gr\"onwall's inequality gives
\begin{equation}\label{est-A-high1}
\begin{split}
\|DA(t)\|_{L^2}\lesssim& \int_0^t \sum_{\alpha'=0}^1 \|u(\tau)-u_0\|_{H^{\alpha'}}\|\nabla \bar a\|_{H^{2-\alpha'}} e^{\int_\tau^t \|u(\tau')\|_{H^2}\, d\tau'}\, d\tau\\
%\lesssim&\ t_N\sum_{\alpha'=0}^1 \|u(t)-u_0\|_{H^{\alpha'}}\|\nabla \bar a\|_{H^{2-\alpha'}}\
\end{split}
\end{equation}
and for $\alpha\geq 2$
\begin{equation}\label{est-A-high2}
\begin{split}
\|D^\alpha A(t)\|_{L^2}\lesssim& \int_0^t\left(\sum_{\alpha'=2}^\alpha \|\nabla u(\tau)\|_{H^{\alpha'}}\|\nabla A(\tau)\|_{H^{\alpha-\alpha'}}\right.\\
&\left.+\sum_{\alpha'=0}^\alpha \|u(\tau)-u_0\|_{H^{\alpha'}}\|\nabla \bar a\|_{H^{\alpha-\alpha'+1}}\right) e^{\int_\tau^t \|u(\tau')\|_{H^2}\, d\tau'}\, d\tau.
%\lesssim&\ t_N\left(\sum_{\alpha'=2}^\alpha \|\nabla u(t)\|_{H^{\alpha'}}\|\nabla A(t)\|_{H^{\alpha-\alpha'}}\right.\\
%&\left.+\sum_{\alpha'=0}^\alpha \|u(t)-u_0\|_{H^{\alpha'}}\|\nabla \bar a\|_{H^{\alpha-\alpha'+1}}\right) \\
\end{split}
\end{equation}

Thanks to Lemma \ref{le-a-osc}, we deduce from \eqref{u-H} that for $0\leq \tau<t\leq t_0$
\begin{equation}\label{u-H2}
\begin{split}
\int_\tau^t \|u(\tau')\|_{H^2}\, d\tau'&\leq C\delta^2t_0^2\lambda^{\frac45(1-\beta)+\frac65\left(6-\beta \right)}\\
&\leq C\delta^2t_0^2 \lambda^{8-2\beta}\\
&\leq CM^{-2}
\end{split}
\end{equation}
where we used \eqref{choice-t-2} again in the last step. Therefore, 
\begin{equation}\label{e-small}
e^{\int_\tau^t \|u(\tau')\|_{H^2}\, d\tau'}\leq e^{CM^{-2}}<2, \ \ \ 0\leq \tau<t\leq t_0
\end{equation}
for large enough $M>0$. Combining \eqref{est-A-high1}, \eqref{est-A-high2} and \eqref{e-small} yields
\begin{equation}\label{est-A-high3}
\begin{split}
\|DA(t)\|_{L^2}\lesssim&\ t_0\sum_{\alpha'=0}^1 \|u(t)-u_0\|_{H^{\alpha'}}\|\nabla \bar a\|_{H^{2-\alpha'}},\\
\|D^\alpha A(t)\|_{L^2}
\lesssim&\ t_0\left(\sum_{\alpha'=2}^\alpha \|\nabla u(t)\|_{H^{\alpha'}}\|\nabla A(t)\|_{H^{\alpha-\alpha'}}\right.\\
&\left.+\sum_{\alpha'=0}^\alpha \|u(t)-u_0\|_{H^{\alpha'}}\|\nabla \bar a\|_{H^{\alpha-\alpha'+1}}\right), \ \ \alpha\geq 2.
\end{split}
\end{equation}

In view of \eqref{est-A-high3} and \eqref{A-l2-final}, we can iterate the process a few times to obtain estimates for $\|DA(t)\|_{L^2}$, $\|D^2A(t)\|_{L^2}$, $\|D^3A(t)\|_{L^2}$ and $\|D^4A(t)\|_{L^2}$. 
%Observing that there are multiple terms on the right hand side in the inequalities of \eqref{est-A-high3}, we first identity the dominating term among them. 
Applying \eqref{u-H} and Lemma \ref{le-a-osc} to \eqref{est-A-high3} gives 
\begin{equation}\label{A-H1-final}
\begin{split}
\|DA(t)\|_{L^2}&\lesssim t_0  \|u(t)-u_0\|_{L^2}\|\nabla \bar a\|_{H^{2}}+t_0  \|u(t)-u_0\|_{H^{1}}\|\nabla \bar a\|_{H^{1}}\\
&\lesssim C_0t_0^2(2M)^{\frac43}\delta^2\lambda^{\frac23(1-\beta)+\frac43(4-\beta)}\delta\left(\delta\lambda^{4-\beta}t_0\right)^3\lambda^{3-\beta}\\
&\quad+C_0t_0^2(2M)^{\frac54}\delta^2\lambda^{\frac34(1-\beta)+\frac54(5-\beta)}\delta\left(\delta\lambda^{4-\beta}t_0\right)^2\lambda^{2-\beta}\\
&\lesssim C_0(2M)^{\frac43} \left(\delta\lambda^{4-\beta}t_0\right)^5\delta\lambda^{1-\beta}+C_0(2M)^{\frac54} \left(\delta\lambda^{4-\beta}t_0\right)^4\delta\lambda^{1-\beta}\\
&\lesssim M^{-2}\delta\lambda^{1-\beta} 
\end{split}
\end{equation}
where we applied \eqref{choice-t-2} in the last step, which indicates $\delta\lambda^{4-\beta}t_0\lesssim M^{-1}$.

Writing the second inequality of \eqref{est-A-high3} with $\alpha=2$ explicitly, we have for $0<t\leq t_0$
\begin{equation}\notag
\begin{split}
\|D^2 A(t)\|_{L^2}
\lesssim&\ t_0 \|\nabla u(t)-\nabla u_0\|_{H^{2}}\|\nabla A(t)\|_{L^{2}}+t_0 \|\nabla u_0\|_{H^{2}}\|\nabla A(t)\|_{L^{2}}\\
&+t_0\sum_{\alpha'=0}^2 \|u(t)-u_0\|_{H^{\alpha'}}\|\nabla \bar a\|_{H^{3-\alpha'}}\\
\lesssim&\ t_0 \|\nabla u(t)\|_{H^{2}}\|\nabla A(t)\|_{L^{2}}+t_0 \|u(t)-u_0\|_{L^2}\|\nabla \bar a\|_{H^{3}}\\
&\quad+...+t_0\|u(t)-u_0\|_{H^2}\|\nabla \bar a\|_{H^{1}}.
\end{split}
\end{equation}
Applying \eqref{A-H1-final} along with Lemma \ref{le-initial}, \eqref{u-H} and Lemma \ref{le-a-osc} to the previous inequality yields 
\begin{equation}\label{A-H2-final}
\begin{split}
\|D^2 A(t)\|_{L^2}&\lesssim C_0t_0^2(2M)^{\frac{7}{6}} \delta^2\lambda^{\frac56(1-\beta)+\frac76(7-\beta)} M^{-2}\delta \lambda^{1-\beta}\\
&\quad+t_0\delta \lambda^{5-\beta}M^{-2}\delta \lambda^{1-\beta}\\
&\quad+C_0t_0^2(2M)^{\frac43}\delta^2\lambda^{\frac23(1-\beta)+\frac43(4-\beta)} \delta(\delta\lambda^{4-\beta}t_0)^4\lambda^{4-\beta}\\
&\quad+C_0t_0^2(2M)^{\frac54}\delta^2\lambda^{\frac34(1-\beta)+\frac54(5-\beta)} \delta(\delta\lambda^{4-\beta}t_0)^3\lambda^{3-\beta}\\
&\quad+C_0t_0^2(2M)^{\frac65}\delta^2\lambda^{\frac45(1-\beta)+\frac65(6-\beta)} \delta(\delta\lambda^{4-\beta}t_0)^2\lambda^{2-\beta}\\
&\lesssim C_0(2M)^{\frac{7}{6}} M^{-2}(\delta \lambda^{4-\beta}t_0)^2 \delta\lambda^{2-\beta}+M^{-2} (\delta \lambda^{4-\beta}t_0)\delta\lambda^{2-\beta}\\
&\quad+ C_0(2M)^{\frac{4}{3}} (\delta \lambda^{4-\beta}t_0)^6 \delta\lambda^{2-\beta}+C_0(2M)^{\frac{5}{4}} (\delta \lambda^{4-\beta}t_0)^5 \delta\lambda^{2-\beta}\\
&\quad +C_0(2M)^{\frac{6}{5}} (\delta \lambda^{4-\beta}t_0)^4 \delta\lambda^{2-\beta}\\
&\lesssim M^{-2}  \delta\lambda^{2-\beta}
\end{split}
\end{equation}
where we used \eqref{choice-t-2} again in the last step.

We iterate to estimate $\|D^3 A(t)\|_{L^2}$ analogously, 
\begin{equation}\notag
\begin{split}
\|D^3 A(t)\|_{L^2}
\lesssim&\ t_0 \|\nabla u(t)\|_{H^2}\|\nabla A(t)\|_{H^1}+t_0\|\nabla u(t)\|_{H^3}\|\nabla A(t)\|_{L^2}\\
&+t_N\sum_{\alpha'=0}^3 \|u(t)-u_0\|_{H^{\alpha'}}\|\nabla \bar a\|_{H^{4-\alpha'}}\\
\lesssim&\ t_0 \|\nabla (u(t)-u_0)\|_{H^2}\|\nabla A(t)\|_{H^1}+t_0\|\nabla u_0\|_{H^2}\|\nabla A(t)\|_{H^1}\\
&\quad+t_0 \|\nabla (u(t)-u_0)\|_{H^3}\|\nabla A(t)\|_{L^2}+t_0 \|\nabla u_0\|_{H^3}\|\nabla A(t)\|_{L^2}\\
&+t_0\sum_{\alpha'=0}^3 \|u(t)-u_0\|_{H^{\alpha'}}\|\nabla \bar a\|_{H^{4-\alpha'}}.
\end{split}
\end{equation}
According to \eqref{A-H1-final} and \eqref{A-H2-final} together with Lemma \ref{le-initial}, \eqref{u-H} and Lemma \ref{le-a-osc}, we continue to infer
\begin{equation}\label{A-H3-final}
\begin{split}
\|D^3 A(t)\|_{L^2}
\lesssim&\ C_0t_0^2(2M)^{\frac76}\delta^2 \lambda^{\frac56(1-\beta)+\frac76(7-\beta)} M^{-2}\delta\lambda^{2-\beta}+t_0\delta \lambda^{5-\beta}M^{-2}\delta\lambda^{2-\beta}\\
&\quad+C_0t_0^2(2M)^{\frac87}\delta^2 \lambda^{\frac67(1-\beta)+\frac87(8-\beta)} M^{-2}\delta\lambda^{1-\beta}+t_0\delta \lambda^{6-\beta}M^{-2}\delta\lambda^{1-\beta}\\
&\quad+C_0t_0^2(2M)^{\frac43}\delta^2\lambda^{\frac23(1-\beta)+\frac43(4-\beta)} \delta(\delta\lambda^{4-\beta}t_0)^5\lambda^{5-\beta}\\
&\quad+C_0t_0^2(2M)^{\frac54}\delta^2\lambda^{\frac34(1-\beta)+\frac54(5-\beta)} \delta(\delta\lambda^{4-\beta}t_0)^4\lambda^{4-\beta}\\
&\quad+C_0t_0^2(2M)^{\frac65}\delta^2\lambda^{\frac45(1-\beta)+\frac65(6-\beta)} \delta(\delta\lambda^{4-\beta}t_0)^3\lambda^{3-\beta}\\
&\quad+C_0t_0^2(2M)^{\frac76}\delta^2\lambda^{\frac56(1-\beta)+\frac76(7-\beta)} \delta(\delta\lambda^{4-\beta}t_0)^2\lambda^{2-\beta}\\
\lesssim&\ M^{-2}\delta \lambda^{3-\beta}
\end{split}
\end{equation}
where we used $\delta\lambda^{4-\beta}t_0\lesssim M^{-1}$ again in the last step.

In the end we proceed to analyze $\|D^4 A(t)\|_{L^2}$, taking $\alpha=4$ in \eqref{est-A-high3}
\begin{equation}\notag
\begin{split}
\|D^4 A(t)\|_{L^2}
\lesssim&\ t_0 \|\nabla (u(t)-u_0)\|_{H^2}\|\nabla A(t)\|_{H^2}+t_0 \|\nabla u_0\|_{H^2}\|\nabla A(t)\|_{H^2}\\
&\quad+t_0 \|\nabla (u(t)-u_0)\|_{H^3}\|\nabla A(t)\|_{H^1}+t_0 \|\nabla u_0\|_{H^3}\|\nabla A(t)\|_{H^1}\\
&\quad+t_0 \|\nabla (u(t)-u_0)\|_{H^4}\|\nabla A(t)\|_{L^2}+t_0 \|\nabla u_0\|_{H^4}\|\nabla A(t)\|_{L^2}\\
&\quad+t_0\sum_{\alpha'=0}^4 \|u(t)-u_0\|_{H^{\alpha'}}\|\nabla \bar a\|_{H^{5-\alpha'}}.
\end{split}
\end{equation}
Applying \eqref{A-H1-final}, \eqref{A-H2-final} and \eqref{A-H3-final} together with Lemma \ref{le-initial}, \eqref{u-H} and Lemma \ref{le-a-osc}, we obtain
\begin{equation}\label{A-H4-final}
\begin{split}
\|D^4 A(t)\|_{L^2}
\lesssim&\ C_0t_0^2(2M)^{\frac76}\delta^2 \lambda^{\frac56(1-\beta)+\frac76(7-\beta)} M^{-2}\delta\lambda^{3-\beta}+t_0\delta \lambda^{5-\beta}M^{-2}\delta\lambda^{3-\beta}\\
&\quad+C_0t_0^2(2M)^{\frac87}\delta^2 \lambda^{\frac67(1-\beta)+\frac87(8-\beta)} M^{-2}\delta\lambda^{2-\beta}+t_0\delta \lambda^{6-\beta}M^{-2}\delta\lambda^{2-\beta}\\
&\quad+C_0t_0^2(2M)^{\frac98}\delta^2 \lambda^{\frac78(1-\beta)+\frac98(9-\beta)} M^{-2}\delta\lambda^{1-\beta}+t_0\delta \lambda^{7-\beta}M^{-2}\delta\lambda^{1-\beta}\\
&\quad+C_0t_0^2(2M)^{\frac43}\delta^2\lambda^{\frac23(1-\beta)+\frac43(4-\beta)} \delta(\delta\lambda^{4-\beta}t_0)^6\lambda^{6-\beta}\\
&\quad+C_0t_0^2(2M)^{\frac54}\delta^2\lambda^{\frac34(1-\beta)+\frac54(5-\beta)} \delta(\delta\lambda^{4-\beta}t_0)^5\lambda^{5-\beta}\\
&\quad+C_0t_0^2(2M)^{\frac65}\delta^2\lambda^{\frac45(1-\beta)+\frac65(6-\beta)} \delta(\delta\lambda^{4-\beta}t_0)^4\lambda^{4-\beta}\\
&\quad+C_0t_0^2(2M)^{\frac76}\delta^2\lambda^{\frac56(1-\beta)+\frac76(7-\beta)} \delta(\delta\lambda^{4-\beta}t_0)^3\lambda^{3-\beta}\\
&\quad+C_0t_0^2(2M)^{\frac87}\delta^2\lambda^{\frac67(1-\beta)+\frac87(8-\beta)} \delta(\delta\lambda^{4-\beta}t_0)^2\lambda^{2-\beta}\\
\lesssim&\ M^{-2}\delta \lambda^{4-\beta}.
\end{split}
\end{equation}
Then, an interpolation of the estimates \eqref{A-H1-final} and \eqref{A-H4-final} yields 
\begin{equation}\label{A-Hs}
\|A(t)\|_{H^s}\leq C_M \delta \lambda^{s-\beta} \ \ \ \forall \ \ t\in[0,t_0]
\end{equation}
for $1<\beta<4$, where $C_M$ is a small constant depending on $M$.

Analogously as in the proof of Lemma \ref{le-high}, we can bootstrap the estimates to obtain
\begin{equation}\label{A-H4-tN}
\|D^s A(t)\|_{L^2}
\leq 2C_M \delta \lambda^{s-\beta}, \ \ t\in[0,t_N].
\end{equation}

%We remark that the condition \eqref{zeta} together with $\gamma>1$ and $0<\zeta<5-\beta$ implies $3<\beta<4$ and $\gamma<\frac{5-\beta}{2(4-\beta)}$.

Regarding $u(t)-u_0$ in $H^{s-2}$, it follows from \eqref{u-u0-basic}, Lemma \ref{le-a-osc} and \eqref{A-Hs} that for $t\in[0,t_N]$
\begin{equation}\notag
\begin{split}
&\quad \|u(t)-u_0\|_{H^{s-2}}\\
&\lesssim \int_0^t\left(\|D\bar a(\tau)\|_{L^2}+\|D A(\tau)\|_{L^2}\right)^{\frac{s}{s+1}}\left(\|D \bar a(\tau)\|_{H^{s+1}}+\|D A(\tau)\|_{H^{s+1}}\right)^{\frac{s+2}{s+1}} d\tau\\
&\lesssim t_N\left( \delta(\delta \lambda^{4-\beta}t_N)\lambda^{1-\beta}+\delta \lambda^{1-\beta}\right)^{\frac{s}{s+1}}
\left( \delta(\delta \lambda^{4-\beta}t_N)^{s+2}\lambda^{s+2-\beta}+\delta \lambda^{s+2-\beta}\right)^{\frac{s+2}{s+1}}\\
&\lesssim t_N\delta^{2-\frac2\beta(s+4)}\lambda^{s+4-2\beta}\\
&\lesssim \delta^{1-\frac{2s+10}\beta}\lambda^{s-\beta}.
\end{split}
\end{equation}
%Straightforward computation shows that for $3<\beta<4$, 
%\[\frac{\beta+1}{\beta^2+5\beta+5}<\frac{10}{41}, \ \ \frac{(\beta+2)^2}{\beta^2+5\beta+5}<\frac{36}{41}.\]
%Combining \eqref{zeta}, we infer
%\begin{equation}\notag
%\zeta>\frac{(\beta+1)(4-\beta)\gamma}{\beta^2+5\beta+5}+\frac{(\beta+2)^2(5-\beta)}{\beta^2+5\beta+5}
%\end{equation}
%and hence the exponent satisfies
%\[(4-\beta)\gamma+\frac{(\beta+2)^2}{\beta+1}(5-\beta)-\frac{\beta^2+5\beta+5}{\beta+1}\zeta<0.\]
In particular, we have
\[\|u(t)-u_0\|_{H^{\beta-2}}\lesssim \delta^{1-\frac{2\beta+10}\beta}, \ \ t\in[0,t_N].\]
It completes the proof of the lemma.

\cbdu

\medskip

\subsection{Proof of the main result Theorem \ref{thm}}
We choose $\lambda$ sufficiently large depending on $\delta$ such that
\[t_N=\delta^{-\frac2\beta-1}\lambda^{\beta-4}<\delta\]
for any $1<\beta<4$.
It is then clear that Theorem \ref{thm} is an immediate consequence of Lemma \ref{le-norm-inflation}, Lemma \ref{le-bar-u} and Lemma \ref{le-perturbation}.

%\bigskip

%\section*{Acknowledgement}
%The author is grateful for ...

%\Endrefs

\bigskip



\begin{thebibliography}{XX}


%\bibliographystyle{plain}
%\bibliography{NS-stability}
%\references {999}

%\bibitem{ADFL}
%M. Acheritogaray, P. Degond, A. Frouvelle and J-G. Liu.
%\newblock {\em Kinetic formulation and global existence for the Hall-Magnetohydrodynamic system}.
%\newblock Kinetic and Related Models, 4: 901--918, 2011.

\bibitem{BT}
C. Bardos and E. S. Titi.
\newblock {\em Loss of smoothness and energy conserving rough weak solutions for the 3D Euler equations}.
\newblock Discrete Contin. Dyn. Syst. Ser. S, 3(2)185--197, 2010.

%\bibitem{Bha}
%A. Bhattacharjee.
%\newblock {\em Impulsive magnetic reconnection in the Earth's magnetotail and the solar corona}.
%\newblock Ann. Rev. Astron. Astrophys., Vol. 42: 365--384, 2004.

\bibitem{BDS}
J. Birn, J. F. Drake, M. A. Shay, B. N. Rogers, R. E. Denton, M. Hesse, M. Kuznetsova, Z. W. Ma, A. Bhattacharjee, A. Otto and P. L. Pritchett.
\newblock {\em Geospace environmental modeling (GEM) magnetic reconnection challenge}.
\newblock J. Geophys. Res., 106, 3715, 2001.

\bibitem{Bis1}
D. Biskamp.
\newblock {\em Magnetic reconnection in plasmas}.
\newblock Cambridge University Press, 2000.

%\bibitem{Bis2}
%D. Biskamp.
%\newblock {\em Magnetohydrodynamic turbulence}.
%\newblock Cambridge University Press, 2003.

\bibitem{BL1}
J. Bourgain and D. Li.
\newblock {\em Strong ill-posedness of the incompressible Euler equation in borderline Sobolev spaces}.
\newblock Invent. Math., 201(1):97--157, 2015.

%\bibitem{BL2}
%J. Bourgain and D. Li.
%\newblock {\em Strong ill-posedness of the incompressible Euler equation in integer $C^m$ spaces}.
%\newblock Geom. Funct. Anal., 25(1):1--86, 2015.

%\bibitem{BDLIS}
%T. Buckmaster, C. De Lellis, P. Isett, and L. Sz\'ekelyhidi.
%\newblock {\em Anomalous dissipation for $1/5$-H\"older Euler flows}.
%\newblock Ann. of Math., Vol. 182, No. 1: 127-172, 2015.

%\bibitem{BDLS}
%T. Buckmaster, C. De Lellis, and L. Sz\'ekelyhidi.
%\newblock {\em Dissipative Euler flows with Onsager-critical spatial regularity}.
%\newblock Comm. Pure Appl. Math., Vol. 69 No. 9, 16131670, 2016.

%\bibitem{BDLSV}
%T. Buckmaster, C. De Lellis, L. Sz\'ekelyhidi, and V. Vicol.
%\newblock {\em Onsager's conjecture for admissible weak solutions}.
%\newblock Comm. Pure Appl. Math., https://doi.org/10.1002/cpa.21781. 2018.

%\bibitem{BV19}
%T. Buckmaster, and V. Vicol.
%\newblock {\em Convex integration and phenomenologies in turbulence}.
%\newblock arXiv: 1901.09023, 2019.

%\bibitem{BV}
%T. Buckmaster, and V. Vicol.
%\newblock {\em Nonuniqueness of weak solutions to the Navier-Stokes equation}.
%\newblock arXiv: 1709.10033v3, 2017.


%\bibitem{BCV}
%T. Buckmaster, M. Colombo, and V. Vicol.
%\newblock {\em Wild solutions of the Navier-Stokes equations whose singular sets in time have Hausdorff dimension strictly less than 1}.
%\newblock arXiv: 1809.00600, 2018.

%\bibitem{BPS}
%S. V. Bulanov, F. Pegoraro and A. S. Sakharov.
%\newblock {\em Magnetic reconnection in electron magnetohydrodynamics}.
%\newblock Phys. Fluids B, 4, 2499, 1992.

%\bibitem{CRW}
%C. Cao, D. Regmi and J. Wu.
%\newblock {\em The 2D MHD equations with horizontal dissipation and horizontal magnetic diffusion}.
%\newblock J. Differential Equations, Vol. 254: 2661--2681, 2013.

\bibitem{CSZ}
L. Chac\'on, A. Simakov and A. Zocco.
\newblock {\em Steady-state properties of driven magnetic reconnection in 2D electron magnetohydrodynamics}.
\newblock Physical Review Letters, 99, 235001, 2007.

%\bibitem{CDL}
%D. Chae, P. Degond and J-G. Liu.
%\newblock {\em Well-posedness for Hall-magnetohydrodynamics}.
%\newblock Ann. Inst. H. Poincar\'e Anal. Non Lineaire, Vol. 31: 555--565, 2014.

%\bibitem{CL}
%D. Chae and J. Lee.
%\newblock {\em On the blow-up criterion and small data global existence for the Hall-magneto-hydrodynamics}.
%\newblock J. Differential Equations, 256: 3835--3858, 2014.

%\bibitem{CS}
%D. Chae,  and M. Schonbek.
%\newblock {\em On the temporal decay for the Hall-magnetohydrodynamic equations}.
%\newblock J. Differential Equations,  Vol. 255: 3971--3982, 2013.

%\bibitem{CWW}
%D. Chae,  R. Wan and J. Wu.
%\newblock {\em Local well-posedness for the Hall--MHD equations with fractional magnetic diffusion}.
%\newblock arXiv:1404.0486v2, 2014.

%\bibitem{CWeng}
%D. Chae and S. Weng.
%\newblock {\em Singularity formation for the incompressible Hall-MHD equations without resistivity}.
%\newblock Ann. I. H. Poincar\'e-AN, Vol. 33: 1009--1022, 2016.

%\bibitem{CW}
%D. Chae and J. Wolf.
%\newblock {\em On partial regularity for the 3D non-stationary Hall magnetohydrodynamics equations on the plane}.
%\newblock Comm. Math. Phys., Vol. 354: 213--230, 2017.

%\bibitem{CDGG}
%J.-Y. Chemin, B. Desjardins, I. Gallagher and E. Grenier.
%\newblock {\em Fluids with anisotropic viscosity}.
%\newblock  Mod\'el. Math. Anal. Num\'er., Vol. 34: 315--335, 2000.

%\bibitem{CD-nse-modes}
%A. Cheskidov and M. Dai.
%\newblock {\em Kolmogorov's dissipation number and the number of degrees of freedom for the 3D Navier-Stokes equations}.
%\newblock  Proceedings of the Royal Society of Edinburg, Section A, Vol. 149, Issue 2: 429--446, 2019.


\bibitem{CD-continuity}
A. Cheskidov and M. Dai.
\newblock {\em Discontinuity of weak solutions to the 3D NSE and MHD equations in critical and supercritical spaces}.
\newblock  Journal of Mathematical Analysis and Applications 481 (2), 123493, 2020.

\bibitem{CD-mhd}
A. Cheskidov and M. Dai.
\newblock {\em Norm inflation for generalized magneto-hydrodynamic system}.
\newblock  Nonlinearity 28 (1), 129, 2014.

\bibitem{CD-norm}
A. Cheskidov and M. Dai.
\newblock {\em Norm inflation for generalized Navier-Stokes equations}.
\newblock  Indiana University Mathematics Journal, 869-884, 2014.



%\bibitem{CDK}
%A. Cheskidov, M. Dai and L. Kavlie.
%\newblock {\em Determining modes for the 3D Navier-Stokes equations}.
%\newblock  Physica D: Nonlinear Phenomena, Vol.374-375:1--9, 2018.

\bibitem{CMZO}
D. C\'ordoba, L. Mart\'inez-Zoroa and W. Oza\'nski.
\newblock {\em Instantaneous gap loss of Sobolev regularity for the 2D incompressible Euler equations}.
\newblock  arXiv: 2210.17458, 2022.

%\bibitem{Dai22}
%M. Dai.
%\newblock {\em Almost sure existence of global weak solutions for supercritical electron MHD}.
%\newblock  arXiv: 2201.08161, 2022.

\bibitem{Dai-emhd-2d}
M. Dai.
\newblock {\em Global existence of 2D electron MHD near a steady state}.
\newblock  arXiv: 2306.13036, 2023.

%\bibitem{Dai1}
%M. Dai.
%\newblock {\em Local well-posedness of the Hall-MHD system in $H^s(\mathbb R^n)$ with $s>\frac n2$}.
%\newblock  Mathemaische Nachrichten. DOI: 10.1002/mana.201800107, 2019.

%\bibitem{Dai2}
%M. Dai.
%\newblock {\em Local well-posedness for the Hall-MHD system in optimal Sobolev spaces}.
%\newblock arXiv: 1803.09556, 2018.

%\bibitem{Dai3}
%M. Dai.
%\newblock {\em Propagation of regularity for the MHD system in optimal Sobolev space}.
%\newblock arXiv: 1707.07754, 2017.

%\bibitem{Dai18}
%M. Dai.
%\newblock {\em Non-unique weak solutions in Leray-Hopf class of the 3D Hall-MHD system}.
%\newblock arXiv: 1812.11311, 2018.

%\bibitem{Dai-hmhd-reg}
%M. Dai.
%\newblock {\em Regularity criterion for the 3D Hall-magneto-hydrodynamics}.
%\newblock Journal of Differential Equations. Vol. 261: 573--591, 2016.

%\bibitem{D}
%M. Dai.
%\newblock {\em Regularity criterion and energy conservation for the supercritical Quasi-Geostrophic equation}.
%\newblock Journal of Mathematical Fluid Mechanics. To appear. ArXiv:1505.02293, 2015.

%\bibitem{DL}
%M. Dai and H. Liu.
%\newblock {\em Long time behavior of solutions to the 3D Hall-magneto-hydrodynamics system with one diffusion}.
%\newblock Journal of Differential Equations, Vol. 266: 7658--7677, 2019.

%\bibitem{DL-well}
%M. Dai and H. Liu.
%\newblock {\em On well-posedness of generalized Hall-magneto-hydrodynamics}.
%\newblock arXiv: 1906.02284, 2019.

\bibitem{Dai-Oh}
M. Dai and S.J. Oh.
\newblock {\em Beale--Kato--Majda-type continuation criteria for Hall- and electron-magnetohydrodynamics}.
\newblock arXiv:2407.04314, 2024.

\bibitem{Dai-W}
M. Dai and C. Wu.
\newblock {\em Dissipation wavenumber and regularity for electron magnetohydrodynamics}.
\newblock Journal of Differential Equations, Vol. 376: 655--681, 2023.


%\bibitem{DSz}
%S. Daneri, and L. Sz\'ekelyhidi.
%\newblock {\em Non-uniqueness and h-principle for H\"older-continuous weak solutions of the Euler equations}.
%\newblock Arch. Ration. Mech. Anal., Vol. 224 No.2: 471--514, 2017.

%\bibitem{DLS1}
%C. De Lellis, and L. Sz\'ekelyhidi.
%\newblock {\em Dissipative continuous Euler flows}.
%\newblock Invent. Math., Vol.193 No. 2: 377--407, 2013.

%\bibitem{DLS2}
%C. De Lellis, and L. Sz\'ekelyhidi.
%\newblock {\em The Euler equations as a differential inclusion}.
%\newblock Ann. of Math.,  Vol.170 No.3: 1417--1436, 2009.


\bibitem{DL}
R. DiPerna and P.L. Lions.
\newblock {\em Ordinary differential equations, transport theory and Sobolev spaces}.
\newblock Invent. Math., 98: 511--547, 1989.

\bibitem{DM}
R. DiPerna and A. J. Majda.
\newblock {\em Oscillations and concentrations in weak solutions of the incompressible fluid equations}.
\newblock Comm. Math. Phys., 108(4): 667--689, 1987.

%\bibitem{DKM}
%J. F. Drake, R. G. Kleva and M. E. Mandt.
%\newblock {\em Structure of thin current layers: implications for magnetic reconnection}.
%\newblock Phys. Rev. Lett. 73, 1251, 1994.

%\bibitem{DS}
%E. Dumas and F. Sueur.
%\newblock {\em On the weak solutions to the Maxwell-Landau-Lifshitz equations and to the Hall-magnetohydrodynamic equations}.
%\newblock Comm. Math. Phys., 330: 1179--1225, 2014.

%\bibitem{EGM}
%T.M. Elgindi, E. Ghoul and N. Masmoudi.
%\newblock {\em Stable self-similar blow-up for a family of nonlocal transport equations}.
%\newblock Analysis and PDE, Vol. 14(3): 891--908, 2021.

\bibitem{EHSX}
T.M. Elgindi, Y. Huang, A.R. Said and C. Xie.
\newblock {\em A classification theorem for steady Euler flows}.
\newblock arXiv: 2408.14662, 2024.

%\bibitem{EJ}
%T.M. Elgindi and I. Jeong.
%\newblock {\em On the effects of advection and vortex stretching}.
%\newblock Archive for Rational Mechanics and Analysis, 235: 1763--1817, 2020.



%\bibitem{FLS}
%D. Faraco, S. Lindberg, and L. Sz\'ekelyhidi.
%\newblock {\em Bounded solutions of ideal MHD with compact support in space-time}.
%\newblock arXiv: 1909.08678, 2019.

%\bibitem{FMRR}
%C. L. Fefferman, D. S. McCormick, J. C. Robinson and J. L. Rodrigo.
%\newblock {\em Higher order commutator estimates and local existence for the non-resistive MHD equations and related models}.
%\newblock Journal of Functional Analysis, Vol. 267: 1035--1056, 2014.

%\bibitem{FHN}
%J. Fan, S. Huang and G. Nakamura.
%\newblock {\em Well-posedness for the axisymmetric incompressible viscous Hall-magnetohydrodynamic equations}.
%\newblock Appl. Math. Lett., 26: 963--967, 2013.



%\bibitem{Fri}
%U. Frisch.
%\newblock {\em Turbulence: The Legacy of A. N. Kolmogrov}.
%\newblock Cambridge University Press, Cambridge, 1995.

%\bibitem{Galtier}
%S. Galtier.
%\newblock {\em Introduction to Modern Magnetohydrodynamics}.
%\newblock Cambridge University Press, London, 2016.


%\bibitem{GH}
%K. N. Gourgouliatos and R. Hollerbach.
%\newblock {\em Resistive tearing instability in electron MHD: application to neutron star crusts}.
%\newblock Monthly Notices of the Royal Astronomical Society, Vol. 463, Iss. 3: 3381--3389, 2016.


%\bibitem{Gr}
%L. Grafakos.
%\newblock {\em Modern Fourier analysis}.
%\newblock Second edition. Graduate Texts in Mathematics, 250. Springer, New York, 2009.



%\bibitem{HouL}
%T. Y. Hou and Z. Lei. 
%\newblock {\em On the stabilizing effect of convection in three-dimensional incompressible flows}.
%\newblock Commun. Pure Appl. Math., 62(4): 501--564, 2009.




%\bibitem{Is0}
%P. Isett.
%\newblock {\em Holder continuous Euler flows with compact support in time}.
%\newblock Ph. D. Thesis, Princeton Jniversity, 2013.

%\bibitem{Is}
%P. Isett.
%\newblock {\em A Proof of Onsager's Conjecture}.
%\newblock Ann. of Math., Vol.188 No.3: 1--93, 2018.


\bibitem{JO1}
I. Jeong and S. Oh.
\newblock {\em On illposedness of the Hall and electron magnetohydrodynamic equations without resistivity on the whole space}.
\newblock arXiv: 2404.13790, 2024.

\bibitem{JO2}
I. Jeong and S. Oh.
\newblock {\em On the Cauchy problem for the Hall and electron magnetohydrodynamic equations without resistivity I: illposedness near degenerate stationary solutions}.
\newblock Annals of PDE, vol.8, no.15, 2022.

\bibitem{JO3}
I. Jeong and S. Oh.
\newblock {\em Wellposedness of the electron MHD without resistivity for large perturbations of the uniform magnetic field}.
\newblock arXiv: 2402.06278, 2024.

%\bibitem{JS1}
%H. Jia and V. \v{S}ver\'ak.
%\newblock {\em Are the incompressible 3d Navier-Stokes equations locally ill-posed in the natural energy space?}
%\newblock J. Funct. Anal., Vol. 268(12): 3734--3766, 2015.

%\bibitem{JS2}
%H. Jia and V. \v{S}ver\'ak.
%\newblock {\em Local-in-space estimates near initial time for weak solutions of the Navier-Stokes equations and forward self-similar solutions}.
%\newblock Invent. Math., Vol. 196(1): 233--265, 2014.

\bibitem{KC}
D.A. Knoll and L. Chac\'on.
\newblock {\em Coalescence of magnetic islands in the low-resistivity, Hall-MHD regime}.
\newblock Phys. Rev. Lett., 96, 135001, 2006.

%\bibitem{K41}
%A. Kolmogoroff.
%\newblock {\em The local structure of turbulence in incompressible viscous fluid for very large {R}eynold's numbers.}
%\newblock {\em C. R. (Doklady) Acad. Sci. URSS (N.S.)}, 30:301--305, 1941.

%\bibitem{Lan} L. D. Landau.
%\newblock {\em On the problem of turbulence}.
%\newblock Doklady Akademii Nauk SSSR. 44: 339--342, 1944.


%\bibitem{L}
%P. G. Lemari\'e-Rieusset.
%\newblock {\em Recent developments in the Navier-{S}tokes problem}.
%\newblock Chapman and Hall/CRC Research Notes in Mathematics, 431. Chapman  and Hall/CRC, Boca Raton, FL, 2002.

%\bibitem{Lei}
%Z. Lei.
%\newblock {\em On axially symmetric incompressible magnetohydrodynamics in three dimensions}.
%\newblock Journal of Differential Equations, Vol. 259, Issue 7: 3202--3215, 2015.

%\bibitem{LXZ}
%F. Lin, L. Xu and P. Zhang.
%\newblock {\em Global small solutions of 2-D incompressible MHD system}.
%\newblock Journal of Differential Equations, Vol. 259, Issue 10: 5440--5485, 2015.

%\bibitem{LZ}
%F. Lin and P. Zhang.
%\newblock {\em Global small solutions to an MHD-type system: the three-dimensional case}.
%\newblock Comm. Pure Appl. Math., Vol. 67: 531--580, 2014.

%\bibitem{Lions}
%J.-L. Lions.
%\newblock {\em Quelques m\'ethodes de r\'esolution des probl\'emes aux limites non lin\'eaires}.
%\newblock volume 1. Dunod; Gauthier-Villars, Paris, 1969.

%\bibitem{LH}
%G. Luo and T.Y. Hou.
%\newblock {\em Potentially singular solutions of the 3D incompressible Euler equations}.
%\newblock PNAS, 111(36): 12968--12973, 2014.

%\bibitem{LT}
%T. Luo and E. Titi.
%\newblock {\em Non-uniqueness of weak solutions to hyperviscous Navier-Stokes equations-On sharpness of J.-L. Lions exponent}.
%\newblock arXiv:1808.07595, 2018.

%\bibitem{Lu}
%X. Luo.
%\newblock {\em Stationary solutions and nonuniqueness of weak solutions for the Navier-Stokes equations in high dimensions}.
%\newblock arXiv:1807.09318, 2018.



%\bibitem{MB}
%A. J. Majda and A. L. Bertozzi.
%\newblock {\em Vorticity and incompressible flow}.
%\newblock Cambridge Jniversity Press, Cambridge, JK, 2001.


%\bibitem{MH}
%H. Miura and D. Hori.
%\newblock {\em Hall effects on local structure in decaying MHD turbulence}.
%\newblock J. Plasma Fusion Res., 8: 73--76, 2009.

%\bibitem{MS}
%S. Modena, and L. Sz\'ekelyhidi.
%\newblock {\em Non-uniqueness for the transport equation with Sobolev vector fields}.
%\newblock arXiv:1712.03867v3, 2017.

%\bibitem{Og}
%T. Ogawa.
%\newblock {\em Sharp Sobolev inequality of logarithmic type and the limiting regularity condition to the harmonic heat flow}.
%\newblock SIAM J. Math. Anal., 34: 1318--1330, 2003.



%\bibitem{On}
%L. Onsager.
%\newblock {\em Statistical hydrodynamics}.
%\newblock Nuovo Cimento (9), 6(Supplemento, 2(Convegno Internazionale di Meccanica Statistica)):279--287, 1949.

%\bibitem{Pai}
%M. Paicu.
%\newblock {\em \'Equation anisotrope de Navier-Stokes dans des espaces critiques}.
%\newblock Rev. Mat. Iberoam., 21: 179--235, 2005.

%\bibitem{PM}
%M. Polygiannakis and X. Mossas.
%\newblock {\em A review of magneto-vorticity induction in Hall- MHD plasmas}.
%\newblock Plasma Phys. Control \& Fusion, 43: 195--221, 2001.

%\bibitem{RY}
%M. M. Rahman and K. Yamazaki.
%\newblock {\em Remarks on the global regularity issue of the two and a half dimensional Hall-magnetohydrodynamics system}.
%\newblock arXiv:2206.12026, 2022.

%\bibitem{RWXZ}
%X. Ren, J. Wu, Z. Xiang and Z. Zhang.
%\newblock {\em Global existence and decay of smooth solution for the 2-D MHD equations without magnetic diffusion}.
%\newblock Journal of Functional Analysis, 267:503--541, 2014.

%\bibitem{RDD}
%B.N. Rogers, R.E. Denton, J.F. Drake and M.A. Shay.
%\newblock {\em ...}.
%\newblock Phys. Rev. Lett., 87, 195004, 2001.

%\bibitem{SJ}
%D. Shalybkov and V. Jrpin.
%\newblock {\em The Hall effect and the decay of magnetic fields}.
%\newblock Astron. Astrophys., 685--690, 1997.

%\bibitem{SDS}
%M.A. Shay, J.F. Drake, M. Swisdak and B.N. Rogers.
%\newblock {\em ...}.
%\newblock Phys. Plasmas, 11, 2199, 2004.

%\bibitem{Tao}
%T. Tao.
%\newblock {\em Finite time blowup for an averaged three-dimensional Navier-Stokes equation}.
%\newblock J. Amer. Math. Soc., Vol. 29: 601--674, 2016.

%\bibitem{Tem}
%R. Temam.
%\newblock {\em Navier-Stokes Equations: Theory and Numerical Analysis}.
%\newblock AMS Chelsea, Providence, Rhode Island, 2000.

%\bibitem{Wal}
%F. Waleffe.
%\newblock {\em On some dyadic models of the Euler equations}.
%\newblock Proc. Amer. Math. Soc., 134 (10): 2913--2922, 2006.

%\bibitem{Wa}
%M. Wardle.
%\newblock {\em Star formation and the Hall effect}.
%\newblock Astrophys. Space Sci., 292: 317--323, 2004.

%\bibitem{WH}
%C.J. Wareing and R. Hollerbach.
%\newblock {\em Forward and inverse cascades in decaying two-dimensional electron magnetohydrodynamic turbulence}.
%\newblock Physics of Plasmas, 16, 042307, 2009.

%\bibitem{WHL}
%T. S. Wood,  R. Hollerbach and  M. Lyutikov.
%\newblock {\em Density-shear instability in electron magneto-hydrodynamics}.
%\newblock Physics of Plasmas, 21, 052110, 2014.

%\bibitem{Zy}
%A. Zygmund.
%\newblock {\em Trigonometric Series}.
%\newblock Cambridge University Press. Third Edition, 2002.


\end{thebibliography}
\end{document}